\magnification 1200

  \input amssym

  % FONTS  ***************************************

  \font \bbfive     = bbm5
  \font \bbseven    = bbm7
  \font \bbten      = bbm10
  \font \eightbf    = cmbx8
  \font \eighti     = cmmi8 \skewchar \eighti = '177
  \font \eightit    = cmti8
  \font \eightrm    = cmr8
  \font \eightsl    = cmsl8
  \font \eightsy    = cmsy8 \skewchar \eightsy = '60
  \font \eighttt    = cmtt8 \hyphenchar\eighttt = -1

  \font \sixi       = cmmi6 \skewchar \sixi = '177
  \font \sixrm      = cmr6
  \font \sixsy      = cmsy6 \skewchar \sixsy = '60
  \font \tensc      = cmcsc10

  \font \titlefont  = cmbx10
  \scriptfont \bffam    = \bbseven
  \scriptscriptfont \bffam = \bbfive
  \textfont \bffam  = \bbten

  \newskip \ttglue

  \def \eightpoint {\def \rm {\fam0 \eightrm }\relax
  \textfont0= \eightrm
  \scriptfont0 = \sixrm \scriptscriptfont0 = \fiverm
  \textfont1 = \eighti
  \scriptfont1 = \sixi \scriptscriptfont1 = \fivei
  \textfont2 = \eightsy
  \scriptfont2 = \sixsy \scriptscriptfont2 = \fivesy
  \textfont3 = \tenex
  \scriptfont3 = \tenex \scriptscriptfont3 = \tenex
  \def \it {\fam \itfam \eightit }\relax
  \textfont \itfam = \eightit
  \def \sl {\fam \slfam \eightsl }\relax
  \textfont \slfam = \eightsl
  \def \bf {\fam \bffam \eightbf }\relax
  \textfont \bffam = \bbseven
  \scriptfont \bffam = \bbfive
  \scriptscriptfont \bffam = \bbfive
  \def \tt {\fam \ttfam \eighttt }\relax
  \textfont \ttfam = \eighttt
  \tt \ttglue = .5em plus.25em minus.15em
  \normalbaselineskip = 9pt
  \def \MF {{\manual opqr}\-{\manual stuq}}\relax
  \let \sc = \sixrm
  \let \big = \eightbig
  \setbox \strutbox = \hbox {\vrule height7pt depth2pt width0pt}\relax
  \normalbaselines \rm }

  % LABELS CONTROL  ***************************************

  % Determines if macro #1 is defined or not.
  \def \ifundef #1{\expandafter \ifx \csname #1\endcsname \relax }

  % Define \showlabel, \showcitations and/or \showlcit if you wanto to see labels,
  % citations, or local-citations, eg.
  % \def \showlabel{}
  % \def \showcitations{}
  % \def \showlcit{}

  \newcount \secno \secno = 0
  \newcount \stno \stno = 0
  \newcount \eqcntr \eqcntr= 0

  % Prints auxiliary information about labels and citations when requested.
  \def \track #1#2#3{\ifundef{#1}\else \hbox{\sixrm[#2\string #3] }\fi}

  % Advances the sequencial numbering system.
  \def \advseqnumbering {\global \advance \stno by 1 \global \eqcntr=0}

  % Prints the sequencial numbering.  \stno is nor printed if zero.
  \def \current {\number \secno \ifnum \number \stno = 0 \else .\number \stno \fi }

  % Prints error message when advanced citations are wrong, and halts execution.
  \def \laberr #1{\edef\a{******** A label has been defined more than once with
conflicting values.  Please check!!! #1 ********}
  {\hfill\break \hrule \noindent \pilar {20pt}\a\stake{20pt} \hrule}
  \message{\a}}

  % Defines a label and, in case it has already been defined, check that new value is same as old.
  \def \deflabel#1#2{%
    \ifundef {#1}%
      \global \expandafter
      \edef \csname #1\endcsname {#2}%
    \else
      \edef\deflabelaux{\expandafter\csname #1\endcsname}%
      \edef\deflabelbux{#2}%
      \ifx \deflabelaux \deflabelbux \else \laberr{#1 ?=? (\deflabelaux) ?=? (\deflabelbux)} \fi
      \fi
    \track{showlabel}{*}{#1}}

  % Defines label for an equation outside a statement.
  \def \eqmark #1 {\advseqnumbering
    \eqno {(\current)}
    \deflabel{#1}{\current}}

  % Defines label for an equation inside a statement.
  \def \subeqmark #1 {\global \advance\eqcntr by 1
    \edef\subeqmarkaux{\current.\number\eqcntr}
    \eqno {(\subeqmarkaux)}
    \deflabel{#1}{\subeqmarkaux}}

  \def \label #1 {\deflabel{#1}{\current}}
  \def \lcite #1{(#1\track{showlcit}{$\bullet$}{#1})}
  \def \forwardcite #1#2{\deflabel{#1}{#2}\lcite{#2}}% Formerly \fcite.

  % CITATIONS CONTROL ***************************************

  % Permits LaTeX style citations e.g. \cite{\paperByFulano} or \cite[Theorem 2.1]{\paperByFulano}.
  \catcode`\@=11
  \def \c@itrk #1{{\bf #1}\track{showcitations}{\#}{#1}} % Cite and track, no brackets
  \def \c@ite #1{[\c@itrk{#1}]}
  \def \sc@ite [#1]#2{[\c@itrk{#2}\hskip 0.7pt:\hskip 2pt #1]}
  \def \du@lcite {\if \pe@k [\expandafter \sc@ite \else \expandafter \c@ite \fi }
  \def \cite {\futurelet \pe@k \du@lcite }
  \catcode`\@=12

  % BIBLIOGRAPHY  ***************************************

  \newcount \bibno \bibno = 0

  \def \bibitem #1#2#3#4{\smallbreak \item {[#1]} #2, ``#3'', #4.}

  \def \references {\begingroup \bigbreak \eightpoint
    \centerline {\tensc References}
    \nobreak \medskip \frenchspacing }

  % HEADER  ***************************************

  \def \Headlines #1#2{\nopagenumbers
    \headline {\ifnum \pageno = 1 \hfil
    \else \ifodd \pageno \tensc \hfil \lcase {#1} \hfil \folio
    \else \tensc \folio \hfil \lcase {#2} \hfil
    \fi \fi }}

  \long \def \Quote #1\endQuote {\begingroup \leftskip 35pt \rightskip 35pt
\parindent 17pt \eightpoint #1\par \endgroup }
  \long \def \Abstract #1\endAbstract {\bigskip \Quote \noindent #1\endQuote}
  \def \Address #1#2{\bigskip {\tensc #1 \par \it E-mail address: \tt  #2}}
  
  \def \Note #1{\footnote {}{\eightpoint #1}}
  \def \Date #1 {\Note {\it Date: #1.}}

  % CONTROL SEQUENCES

  \def \lcase #1{\edef \auxvar {\lowercase {#1}}\auxvar }

  \def \section #1 \par{\global \advance \secno by 1 \stno = 0
    \bigbreak \noindent {\bf \number \secno .\enspace #1.}
    \nobreak \medskip \noindent}

  % \sl is the default font for \state, but it may be overridden
  \def \state #1 #2\par{\medbreak \noindent \advseqnumbering {\bf \current.\enspace #1.\enspace \sl #2\par }\medbreak }
  \def \definition #1\par {\state Definition \rm #1\par }

  % If a proof ends in a displayed equation, \endProof provides for the closing "$$"
  \long \def \Proof #1\endProof {\medbreak \noindent {\it Proof.\enspace }#1
\ifmmode \eqno \endproofmarker $$ \else \hfill $\endproofmarker$ \looseness = -1 \fi \medbreak}

  \def \$#1{#1 $$$$ #1}

  \def \pilar #1{\vrule height #1 width 0pt}
  \def \stake #1{\vrule depth #1 width 0pt}

  \newcount \footno \footno = 1
  \newcount \halffootno \footno = 1
  \def \footcntr {\global \advance \footno by 1
  \halffootno =\footno
  \divide \halffootno by 2
  $^{\number\halffootno}$}
  \def \fn#1{\footnote{\footcntr}{\eightpoint#1\par}}

  % ITENS  ***************************************

  \def \Item #1{\smallskip \item {{\rm #1}}}
  \newcount \zitemno \zitemno = 0

  \def \izitem {\zitemno = 0}
  \def \zitemplus {\global \advance \zitemno by 1 \relax}
  \def \rzitem{\romannumeral \zitemno}
  \def \rzitemplus {\zitemplus \rzitem}
  \def \zitem {\Item {{\rm(\rzitemplus)}}}
  \def \zitemmark #1 {\deflabel{#1}{\rzitem}}

  \newcount \nitemno \nitemno = 0
  \def \initem {\nitemno = 0}
  \def \nitem {\global \advance \nitemno by 1 \Item {{\rm(\number\nitemno)}}}

  \newcount \aitemno \aitemno = -1
  \def \boxlet#1{\hbox to 6.5pt{\hfill #1\hfill}}
  \def \iaitem {\aitemno = -1}
  \def \aitemconv{\ifcase \aitemno a\or b\or c\or d\or e\or f\or g\or
h\or i\or j\or k\or l\or m\or n\or o\or p\or q\or r\or s\or t\or u\or
v\or w\or x\or y\or z\else zzz\fi}
  \def \aitem {\global \advance \aitemno by 1\Item {(\boxlet \aitemconv)}}
  \def \aitemmark #1 {\deflabel{#1}{\aitemconv}}

  % STANDARD DEFINITIONS

  \font\mf=cmex10
  \def\union {\mathop{\raise 9pt \hbox{\mf S}}\limits}
  \def\inters{\mathop{\raise 9pt \hbox{\mf T}}\limits}

  \def \<{\left \langle \vrule width 0pt depth 0pt height 8pt }
  \def \>{\right \rangle }
  \def \ds{\displaystyle}
  \def \and {\hbox {,\quad and \quad }}
  \def \calcat #1{\,{\vrule height8pt depth4pt}_{\,#1}}
  
  \def \imply {\mathrel{\Rightarrow}}
  \def \for #1{,\quad \forall\,#1}
  \def \endproofmarker {\square}
  \def \"#1{{\it #1}\/}
  \def \inv {^{-1}}
  \def \*{\otimes}
  \def \caldef #1{\global \expandafter \edef \csname #1\endcsname {{\cal #1}}}
  \def \N {{\bf N}}

  % PICTEX ***************************************
  %
  \input pictex

  % RANDOM DEFS ***************************************

  \font \rs = rsfs10

  \caldef R
  \caldef I
  \caldef S
  \def\F{{\bf F}}
  \def\Fmn{{\F_{m+n}}}
  \def\O#1#2{{\cal O}_{#1,#2}}
  \def\Omn{\O mn}
  \def\red{r}
  \def\Omnr{\Omn^\red}
  \def\Lmn{L_{m,n}}
  \def\univ{^{u}}
  \def\OR{\Omega\univ}

  \def\s{\sigma}
  \def\phi{\varphi}
  \def\sep{\kern 20pt}
  \def\vsep{20pt}
  \def\us{\underline{s}}
  \def\ut{\underline{t}}
  \def\curly#1{\hbox {\rs #1\kern2pt}}
  \def\.{\odot}
  \def\B{\curly B}
  \def\C{\curly C}
  \def\Cr#1{C^*_{\red}(#1)}
  \def\clspan{\overline{\rm span}}
  \def\BigB{B\big(H\*K\*\ell_2(G)\big)}
  \def\lto#1{\mathrel{\buildrel #1 \over \to}}
  \def\Lto#1{\mathrel{\buildrel #1 \over \longrightarrow}}
  \def\HighLto#1{\mathrel{\buildrel {\raise 2pt \hbox{$\scriptstyle #1$}} \over \longrightarrow}}

  \def\phi{\varphi}
  \def\Cplx{{\bf C}}
  \def\onto{\to \kern -8pt \to}
  \def\i{\subseteq}

  \def\a{\alpha}
  \def\d{\delta}
  \def\b{\beta}

  \caldef H
  \caldef V

  \def\compos{\,{\circ}\,}

  \def\Lin#1{\curly L\big(\ell_2(#1)\big)}

  \def\l{\lambda}
  \caldef I
  \def\IO{\I(\Omega)}
  \def\tu{\theta\univ}
  \def\xr{\rtimes^{\red}}
  \def\MCP{C(\OR)\rtimes_{\tu}\Fmn}
  \def\CP{C(\Omega)\rtimes_{\theta}\Fmn}
  \def\tm{{\ds \mathop{\*}\limits_{\raise 2pt\hbox{\fiverm max}}}}
  \def\Cf#1{C^*(#1)}
  \def\mnsys{$(m,n)$--dynamical system}
  \def\resfinAlg{residually finite}
  \def\resfinGrp{residually finite-dimensional}

  % REFERENCES ***************************************************************************

  \font\bibfont = cmbx8

  \def\AranAbrams{{\bibfont AA}}
  \def\ArGo{{\bibfont AG1}}
  \def\AraGood{{\bibfont AG2}}
  \def\ArMorPa{{\bibfont AMP}}
  \def\FD{{\bibfont FD}}
  \def\LBrown{{\bibfont B}}
  \def\Nate{{\bibfont BO}}
  \def\CuntzKrieger{{\bibfont CK}}
  \def\Hazrat{{\bibfont H}}
  \def\Leavitt{{\bibfont L}}
  \def\Pat{{\bibfont DP}}
  \def\Cohen{{\bibfont C}}
  \def\ExelCircle{{\bibfont E1}}
  \def\amena{{\bibfont E2}}
  \def\ortho{{\bibfont E3}}
  \def\exact{{\bibfont E4}}
  \def\tpa{{\bibfont E5}}
  \def\infinoa{{\bibfont EL}}
  \def\topfree{{\bibfont ELQ}}
  \def\KR{{\bibfont KR}}
  \def\McClanahanUnitary{{\bibfont M1}}
  \def\McClanahanReactangle{{\bibfont M2}}
  \def\McClanahanAmalgamated{{\bibfont M3}}
  \def\McClanahan{{\bibfont M4}}

  % Start of the paper *******************************************************************

  \Headlines {DYNAMICAL SYSTEMS OF TYPE $(m,n)$}
    {P.~Ara, R.~Exel and T.~Katsura}
  \null\vskip -1cm
  \centerline{\titlefont DYNAMICAL SYSTEMS OF TYPE $(m,n)$ AND}
  \smallskip
  \centerline{\titlefont THEIR C*-ALGEBRAS}

  \Date{20 April 2011}

  \footnote{\null}
  {\eightrm 2010 \eightsl Mathematics Subject Classification:
  \eightrm
  46L05, % General theory of C*-algebras
  46L55. % NC Dyn syst.
  }

  \Note
  {\it Key words and phrases: \rm Leavitt C*-algebra, {\mnsys}, exact
C*-algebra, Fell bundles, partial representation, partial action, crossed
product, free group. }

  \bigskip
  \centerline{\tensc Pere Ara, Ruy Exel and Takeshi Katsura}

  \Note{The first-named author was partially supported by DGI MICIIN-FEDER MTM2008-06201-C02-01, and
by the Comissionat per Universitats i Recerca de la Generalitat de Catalunya.
The second-named author was partially supported by CNPq.  The third-named author was partially supported
by the Japan Society for the Promotion of Science.}

  \bigskip

\Abstract
  Given positive integers $n$ and $m$, we consider dynamical systems in which $n$ copies
of a topological space is homeomorphic to $m$ copies of that same space.
  The universal such system is shown to arise naturally from the study of a
C*-algebra we denote by $\Omn$, which in turn is obtained as a quotient of the
well known Leavitt C*-algebra $\Lmn$, a process meant to transform
the generating set of partial isometries of $\Lmn$ into a tame set.
  Describing $\Omn$ as the crossed-product of the universal {\mnsys} by a partial
action of the free group $\Fmn$, we show that $\Omn$ is not exact when $n$ and
$m$ are both greater than or equal to 2, but the corresponding reduced
crossed-product,  denoted $\Omnr$,  is shown to be exact and non-nuclear.
  Still under the assumption that $m, n\geq2$, we prove that the partial action
of $\Fmn$ is topologically free and that $\Omnr$ satisfies property (SP) (small
projections).
  We also show that  $\Omnr$ admits no
finite dimensional representations.
  The techniques developed to treat this system include several new results
pertaining to the theory of Fell bundles over discrete groups.  \endAbstract

\section Introduction

  The well known one-sided shift on $n$ symbols is a dynamical system in which
the configuration space is homeomorphic to $n$ copies of itself.  In this paper
we study systems in which $n$ copies of a topological space $Y$ is homeomorphic
to $m$ copies of it.

Precisely,  this means that one is given a pair $(X,Y)$ of compact Hausdorff topological spaces
$(X,Y)$ such that
  $$
  X = \bigcup_{i=1}^n H_i = \bigcup_{j=1}^m V_j,
  $$
  where the $H_i$ are pairwise disjoint clopen subsets of $X$, each of which is
homeomorphic to $Y$ via given homeomorphisms
  $
  h_i: Y \to H_i,
  $
  and the $V_i$ are pairwise disjoint clopen subsets of $X$, each of which is
homeomorphic to $Y$ via given homeomorphisms
  $
  v_i: Y \to V_i.
  $

\vskip -10pt
\noindent \hfill
\beginpicture
\setcoordinatesystem units <0.0030truecm, 0.0030truecm> point at 3000 0
\setplotarea x from -1700 to 1000, y from -1400 to 1000
\font\bigmath cmmi12 scaled 1440

\putrectangle corners at -1000 200 and -1400 -200
\put {\bigmath Y} at -1600 000

\plot 250 500 250 -500 / \put {$V_1$} <-6pt,8pt> at 200 500
\plot 500 500 500 -500 / \put {$V_2$} <-6pt,8pt> at 500 500
                         \put {$\cdots$} <-6pt,8pt> at 750 500
\plot 750 500 750 -500 / \put {$V_m$} <-6pt,8pt> at 1000 500

\putrectangle corners at 000 500 and 1000 -500 \put {\bigmath X} at 1400 00
\plot 000 200 1000 200 / \put {$H_1$} <10pt,12pt> at 1000 200
                          \put {$\vdots$} <10pt,12pt> at 1000 -100
\plot 000 -200 1000 -200 / \put {$H_n$} <10pt,12pt> at 1000 -500

\arrow  <0.15cm> [0.25,0.75] from -945 100 to -100 380 \put {$h_1$} at -300 430

\put {$\vdots$} at -300 50

\arrow  <0.15cm> [0.25,0.75] from -945 -100 to -100 -380 \put {$h_n$} at -300 -430

\setquadratic
\plot -980 -250 -400 -800 100 -600 / \arrow  <0.15cm> [0.25,0.75] from 75 -630 to 100 -600
                                     \put {$v_1$} at -130 -700
                                     \put {$\ldots$} at 450 -700
                                     \put {$v_m$} at 900 -800
\plot -1200 -250 -00 -1200 850 -600 / \arrow  <0.15cm> [0.25,0.75] from 835 -630 to 850 -600

\endpicture
\hfill\null

\centerline{\eightrm Diagram
  \advseqnumbering
    (\current)
    \deflabel{DiagramOne}{\current}}

\bigskip
\bigskip

To the quadruple
  $\big(X,Y,\{h_i\}_{i=1}^n,\{v_j\}_{j=1}^m\big)$
  we give the name of an {\mnsys}.  When $n$ or $m$ are 1, this
essentially reduces to the shift, but when $m,n\ge 2$, a very different
behavior takes place.

The origin of the ideas developed in the present paper can be traced
back to the seminal work of Cuntz and Krieger \cite{\CuntzKrieger},
where a dynamical interpretation of the Cuntz-Krieger C*-algebras is
given. In particular, the Cuntz algebra ${\cal O} _n$ corresponds to
the full shift on $n$ symbols. Since we are using an ``external"
model for this dynamical system, the C*-algebra $\O{1}{n}$ that we
attach to the $(1,n)$-dynamical system is isomorphic to the algebra
$M_2({\cal O}_n)$.

From a purely algebraic perspective, a motivation to study such
systems comes from the study of certain rings constructed by Leavitt
\cite{L} with the specific goal of having the free module of rank
$n$ be isomorphic to the free module of rank $m$. We refer the
reader to \cite{\AranAbrams}, \cite{\ArMorPa}, \cite{\ArGo},
\cite{\Hazrat} for various interpretations and generalizations of
the algebras constructed by Leavitt to the setting of graph
algebras.

A similar idea lies behind the investigations conducted by Brown
\cite{\LBrown} and McClanahan \cite{\McClanahanUnitary},
\cite{\McClanahanReactangle}, \cite{\McClanahanAmalgamated}, on the C*-algebras
$U^{\rm nc}_{m,n}$. These are the C*-algebras generated by the
entries of a universal unitary matrix of size $m\times n$. It has
been observed in \cite{\AraGood} that there are isomorphisms
$$L_{m,n}\cong M_{m+1}(U^{\rm nc}_{m,n})\cong M_{n+1} (U^{\rm nc}_{m,n}),$$
where $L_{m,n}$  is the universal C*-algebra generated by partial
isometries
  $$
  s_1,\ldots,s_n, \ t_1,\ldots,t_m,
  $$
  sharing the same source projection, and such that the sum of the
range projections of the $s_i$, as well as that of the $t_j$, add up to the
complement of the common source projection.  Incidentally $\Lmn$
may also be constructed as a separated graph C*-algebra
\cite{\AraGood}.

The partial isometries generating this algebra have a somewhat
stubborn algebraic behavior, not least because their final
projections fail to commute. Sidestepping this very delicate issue
we choose to mod out all of the nontrivial commutators and,
  after performing this perhaps rather drastic transformation, we are left with a
C*-algebra which we denote by $\Omn$, and which is consequently generated by
a \"{tame} (see definition \forwardcite{DefineTame}{2.2} below) set of partial
isometries.

We then take advantage of the existing literature on C*-algebras generated by
tame sets of partial isometries \cite{\topfree, \ortho, \infinoa} to describe
$\Omn$ as the crossed product associated to a partial action $\tu$ of the free
group $\Fmn$ on
a compact space $\OR$.  In symbols
  $$
  \Omn \simeq \MCP.
  $$
  It is perhaps no coincidence that the above partial action of $\Fmn$
is given by
an {\mnsys}, as defined above, which is in fact the universal one
\forwardcite{UniversalPropertyOfSystem}{3.8}.

  While our description of the universal {\mnsys} $\OR$ as a subset of the power
set of $\Fmn$ is satisfactory for some purposes, its tree-like
structure may not make it easy to be studied from some points of
view.  We therefore present an alternative version of it in terms of
functions defined on a certain space of finite paths \forwardcite
{FuncEqualConfig}{4.1}.  With this description at hand we are able
to show that the partial action of $\Fmn$ on $\OR$ is topologically
free \forwardcite{Onmtopfree}{4.6}.  When $1\le m,n$ we show that
every nonzero hereditary subalgebra of $\Omnr$ contains a nonzero
projection belonging to $C(\OR)$.

We then initiate a systematic study of $\Omn$, begining with the
fundamental questions of nuclearity and exactness (see \cite{\Nate}
for an extensive study of these important properties of
C*-algebras).

When either $n=1$, or $m=1$, these algebras are Morita--Rieffel equivalent to
Cuntz algebras, so we concentrate on the case in which $n$ and $m$ are greater
than or equal to 2.  Under this condition we prove that $\Omn$ is not nuclear,
and not even exact \forwardcite{OnmNotExact}{7.2}.  However, when we pass to its
\"{reduced} version, namely the reduced crossed product
\cite{\McClanahan}
  $$
  \Omn^\red = C(\OR)\xr_{\tu}\Fmn
  $$
  we find that $\Omn^\red$ is exact, although still not nuclear.

Since the crossed product by a partial action may be defined as the
cross-sectional C*-algebra of the semidirect product Fell bundle, we dedicate a
significant amount of attention to these and in fact many of our statements
about $\Omn$ or $\Omn^\red$ come straight from corresponding results we prove
for general Fell bundles.

If $\B$ is a Fell bundle over a discrete exact group whose unit fiber is an
exact C*-algebra, we prove in \forwardcite{ExactReducedCrossSec}{5.2} that the
reduced cross-sectional C*-algebra $\Cr\B$ is exact. From this it follows that
the reduced crossed product of an exact C*-algebra by a partial action of an exact group is exact,
and hence that $\Omn^\red$ is exact.

Being $\Omn$ a \"{full} crossed product, we are led to study full
cross-sectional C*-algebras of Fell bundles.  The well known fact
\cite[10.2.8]{\Nate} that the maximal tensor product of the reduced
group C*-algebra by itself contains the full group C*-algebra is
generalized in \forwardcite{CrazyFact}{6.2}, where we prove that if
$\B$ is a Fell bundle over the group $G$, then the full
cross-sectional C*-algebra $\Cf\B$ is a subalgebra of $\Cr\B \tm \Cr
G$.  As an immediate consequence we deduce that, if $\Cr\B$ is
nuclear, then the full and reduced cross-sectional C*-algebras of
$\B$ agree \forwardcite{NuclearAmenable}{6.4}.

As another Corollary of \forwardcite{CrazyFact}{6.2} we prove that, if $H$ is a
subgroup of $G$, then the natural map from the full cross-sectional C*-algebra
of the bundle restricted to $H$ embeds in the algebra for the whole bundle.
This result turns out to be crucial in our proof that, in a
partial action, every {\resfinGrp} isotropy group is amenable when the full
cross-sectional algebra is exact \forwardcite{LotsOfConclusions}{7.1}.

When $m,n\geq2$, we show that there are non-amenable
\forwardcite{OnmNotExact}{7.2} isotropy groups in the universal {\mnsys}, so
exactness of $\Omn$ is ruled out by \forwardcite{LotsOfConclusions}{7.1}.

We also consider the question of existence of finite dimensional representations
of $\Omn$ and of $\Omnr$.  A trivial argument
\forwardcite{EasyInexistFinDimRep}{8.1} proves that, when $n\neq m$, neither
$\Omn$ nor $\Omnr$ admit finite dimensional representations.

The case $m=n$ is however a lot more subtle.  While it is easy to produce many
finite dimensional representations of $\Omn$, we have not been able to decide
whether or not there are enough of these to separate points.  In other
words we have not been able to decide whether $\Omn$ is {\resfinAlg}.

With respect to $\Omnr$, we settle the question in \forwardcite
{NoFinDimRepReg}{9.5}, proving that $\Omnr$ admits no finite dimensional
representation for all $m, n\geq2$.

\section The Leavitt C*-algebra

Throughout this paper we fix positive integers $n$ and $m$, with
$m\le n$.

\definition
  The \"{Leavitt C*-algebra of type $(m,n)$} is the universal unital C*-algebra
$\Lmn$ generated by partial isometries
  $s_1,\ldots,s_n, \ t_1,\ldots,t_m$
  satisfying the relations
  $$
  \left. \matrix{
  s_i^*s_{i'} = 0, \hbox{ for } i\neq i',                   \hfill \pilar{0pt}\cr
  t_j^*t_{j'} = 0, \hbox{ for } j\neq j',                   \hfill \pilar{\vsep}\cr
  \ds s_i^*s_i = t_j^*t_j =: w,                         \hfill \pilar{\vsep}\cr
  \ds\sum_{i=1}^n s_is_i^* = \sum_{j=1}^m t_jt_j^* =:v, \hfill \pilar{\vsep}\cr
  \ds vw=0, \quad v+w=1.                                \hfill \pilar{13pt}}
  \right\}{(\R)}
  $$

By choosing a specific representation, it is not difficult to see
that $s_1s_1^*$ does not commute with $t_1t_1^*$ when $m, n\ge2$,  and hence that
$s_1^*t_1$  is not a partial isometry (see e.g.~\cite[5.3]{\ortho}).
This is in contrast with many well known examples of C*-algebras
generated by sets of partial isometries which are almost always
\"{tame} according to the following:

\definition \label DefineTame
  A set $U$ of partial isometries in a C*-algebra is said to be \"{tame} if
every element of $\langle U\cup U^*\rangle$ (meaning the multiplicative
semigroup generated by $U\cup U^*$) is a partial isometry.

See \cite[5.4]{\ortho} for equivalent conditions characterizing tame sets of
partial isometries.

The standard partial isometries generating the Cuntz--Krieger algebras
form a tame set \cite[5.2]{\amena}, as do the corresponding ones for graph
C*-algebras, higher rank graph C*-algebras and many others.

Rather than attempt to face the wild set of partial isometries in $\Lmn$
(incidentally a task not everyone shies away from \cite{\AraGood}), we will
force it to become tame by considering a quotient of $\Lmn$.
  In what follows we will denote by $U$ the subset of partial isometries in
$\Lmn$ that is most relevant to us,  namely
  $$
  U = \{  s_1,\ldots,s_n, \ t_1,\ldots,t_m\}.
  $$

\definition  \label DefineOnm
  We will let $\Omn$ be the quotient of $\Lmn$ by the closed two-sided ideal
generated by all elements of the form
  $$
  xx^*x - x,
  $$
  as $x$ runs in $\langle U\cup U^*\rangle$.  We will denote the images of the
$s_i$ and the $t_j$ in $\Omn$ by $\us_i$ and $\ut_j$, respectively.

It is therefore evident that
  $$
  \{\us_1,\ldots,\us_n, \ \ut_1,\ldots,\ut_m\}
  $$
  is a tame set of partial isometries. In fact it is not hard to prove that
$\Omn$ is the universal unital C*-algebra generated by a tame set of partial
isometries
  satisfying relations $(\R)$.

Let $\Fmn$ denote the free group generated by a set with  $m+n$ elements, say
  $$
  \{a_1,\ldots,a_n, \ b_1,\ldots,b_m\}.
  $$
  Using \cite[5.4]{\ortho} we conclude that there exists a (necessarily unique)
semi-saturated \cite[5.3]{\ortho} partial representation
  $$
  \s :\Fmn \to \Omn
  $$
  such that
  $\s(a_i) = \us_i$, and $\s(b_j) = \ut_j$ (when stating conditions such as
these, which are supposed to hold for every $i=1, \ldots, n$,  and every $j=1,
\ldots, m$,  we will omit making explicit reference of the range of variation of
$i$ and $j$,  which should always be understood as being $1$--$n$,  and $1$--$m$, as above).

Another universal property enjoyed by $\Omn$ is described next.

\state Proposition \label FirstUnivResult
  Let $\rho$ be a semi-saturated partial representation of\/ $\Fmn$ in a unital
C*-algebra $B$ such that the elements
  $s'_i:= \rho(a_i)$ and $t'_j:= \rho(b_j)$ satisfy relations $(\R)$.  Then
there exists a unique unital *-homomorphism $\phi: \Omn \to B$ such that $\rho =
\phi\circ\s$.

\Proof
  Since $\rho$ is a partial representation, one has that the $s'_i$ and the
$t'_j$ are partial isometries.
   By universality of $\Lmn$ one concludes that there exists a unital
*-homomorphism $\psi:\Lmn \to B$, such that
  $\psi(s_i) = s'_i$, and $\psi(t_j) = t'_j$.

Observe that if $x$ is in $\langle U\cup U^*\rangle$, then $\psi(x)$ lies in the
multiplicative semigroup generated by the $s'_i$, the $t'_j$, and their
adjoints.  Employing \cite[5.4]{\ortho} we have that $\psi(x)$ is a partial
isometry and hence that
  $\psi(xx^*x-x)=0$.
  This implies that $\psi$ vanishes on the ideal referred to in
\lcite{\DefineOnm} and hence that it factors through $\Omn$ providing a
*-homomorphism $\phi:\Omn \to B$, such that
  $\phi(\us_i)=s'_i$, and $\phi(\ut_j)=t'_j$.  Therefore
  $$
  \phi(\s(a_i))=   \phi(\us_i)= s'_i = \rho(a_i),
  $$
  and similarly $\phi(\s(b_j)) = \rho(b_j)$.  In other words, $\phi\circ\s$
coincides with $\rho$ on the generators of $\Fmn$.  Since both $\s$ and $\rho$
are semi-saturated, we now conclude that $\phi\circ\s= \rho$ on the whole of
$\Fmn$.  \endProof

So $\Omn$ is the universal unital C*-algebra for partial representations of $\Fmn$
subject to the relations $(\R)$, according to \cite[4.3]{\topfree}, and hence we
may apply \cite[4.4]{\topfree} to deduce that there exists a certain partial
dynamical system
  $(\OR,\Fmn,\tu)$ and
  a *-isomorphism
  $$
  \Psi: \Omn \to \MCP.
  \eqmark IsoOnmProd
  $$
  The choice of notation, specifically the use of the superscript ``$u$",
is motivated by \"{universal} properties to be described below.
  Before giving further details on the above result let us introduce a variation
of $\Omn$.

\definition \label DefineRedOnm
  For every pair of positive integers $(m,n)$ we shall let $\Omn^\red$
denote the corresponding \"{reduced} crossed product
  $$
  \Omn^\red = C(\OR)\xr_{\tu}\Fmn.
  $$

For the convenience of the reader we will now give a brief description of $\OR$
and of the partial action $\tu$.  We refer the reader to \cite[Section
4]{\topfree} for further details.

The first step is to write the relations
defining our algebra in terms of the \"{final projections}
  $$
  e(g) := \s(g)\s(g\inv),
  $$
  for $g\in\Fmn$.  Once this is done we arrive at
  $$
  \left. \matrix{
  e(a_i)e(a_{i'})=0, \hbox{ for } i\neq i',                   \hfill \pilar{0pt}\cr
  e(b_j)e(b_{j'})=0, \hbox{ for } j\neq j',                   \hfill \pilar{\vsep}\cr
  e(a_i\inv) = e(b_j\inv) =: w,                           \hfill \pilar{\vsep}\cr
  \ds\sum_{i=1}^n e(a_i) = \sum_{j=1}^m e(b_j) =:v,       \hfill \pilar{\vsep}\cr
  \ds vw=0, \quad v+w=1.                                  \hfill \pilar{13pt}}
  \right\}{(\R')}
  $$

Observe that, since the $e(g)$ are projections, all of the
above relations expressing orthogonality,  that is,  those having a zero as the
right-hand-side,  follow from ``$v+w=1$".

If we are to apply the theory of \cite[Section 4]{\topfree} to our algebra, we
need to add another relation to $(\R')$ in order to account for the fact that
the partial representations involved in \lcite{\FirstUnivResult} are required to
be semi-saturated.  Although the definition of semi-saturatedness, namely
  $$
  |hk| = |h|+|k| \ \imply \ \s(hk)=\s(h)\s(k),
  $$
  is not expressed in terms of the $e(g)$, we may use \cite[5.4]{\amena} to
replace it with the equivalent form
  $$
  |hk| = |h|+|k| \ \imply \ e(hk)\leq e(h).
  $$

The next step is to translate each of the above relations in terms of equations
on $\{0,1\}^{\Fmn}$.  For this we will find it convenient to identify this product
space with the power set $\curly{P}(\Fmn)$ in the usual way.

According to \cite[Section 4]{\topfree} and \cite[Section 2]{\infinoa} the
translation process consists in replacing each occurrence of a final projection
$e(g)$ in the above relations with the scalar valued function $1_g$ defined by
  $$
  1_g: \xi \in \{0,1\}^{\Fmn} \mapsto [g\in\xi].
  $$
  Here we use brackets to denote Boolean value and we see the truth values ``1" and
``0" as complex numbers.  Therefore $1_g$ is nothing but the characteristic
function of the set
  $$
  \big\{\xi\in \{0,1\}^{\Fmn}: g\in\xi\big\}.
  $$

The description of $\OR$ given
in \cite[4.1]{\topfree} therefore becomes: a necessary and sufficient condition for a
given $\xi\in \{0,1\}^\Fmn$ to belong to $\OR$ is that
  $
  1\in\xi,
  $
  and that
  $$
  \def\gxi{{(g\inv\xi)}}
  \def\vsep{15pt}
  \left. \matrix{
  \left( 1_{hk}1_h-1_{hk}\right)\gxi =0, \hbox{ whenever }  |hk| = |h|+|k|,   \hfill \pilar{0pt}\cr
  \left( 1_{a_i}1_{a_{i'}}\right)\gxi=0, \hbox{ for } i\neq i',                \hfill \pilar{\vsep}\cr
  \big( 1_{b_j}1_{b_{j'}}\big)\gxi=0, \hbox{ for } j\neq j',                   \hfill \pilar{\vsep}\cr
  1_{a_i\inv}\gxi = 1_{b_j\inv}\gxi =: w\gxi,                                  \hfill \pilar{\vsep}\cr
  \ds\sum_{i=1}^n 1_{a_i}\gxi = \sum_{j=1}^m 1_{b_j}\gxi =:v\gxi,       \hfill \pilar{\vsep}\cr
  \ds \left( vw\right)\gxi=0, \quad \left( v+w\right)\gxi=1,                   \hfill \pilar{10pt}}
  \right\}{(\R'')}
  $$
  for every $g\in\xi$.

  For example, to account for the second equation above, it is required that
  $$
  0 =  (1_{a_i}1_{a_{i'}})(g\inv\xi) =
  1_{a_i}(g\inv\xi)\ 1_{a_{i'}}(g\inv\xi) =
  \big[a_i\in g\inv\xi\big]\big[a_{i'}\in g\inv\xi\big] \$=
  \big[ga_i\in \xi\big]\big[ga_{i'}\in \xi\big] =
  \big[ga_i\in \xi\  \wedge \ ga_{i'}\in \xi\big].
  $$
  This may be interpreted as saying that for every $g\in\xi$, not more than one
element of the form $ga_i$ belongs to $\xi$.

As another example, recall from that \cite[4.5]{\infinoa} that, in order for
$\xi$ to satisfy the conditions related to the first equation in $(\R'')$, it is
required that $\xi$ be convex \cite[4.4]{\infinoa}.

The reader may now check that the elements of $\OR$ are precisely
those $\xi\subseteq \Fmn$ such that
  \iaitem

  \aitem $1\in\xi$,

  \aitem $\xi$ is convex,

  \aitem for any $g\in\xi$, one and only one of the conditions below are satisfied:

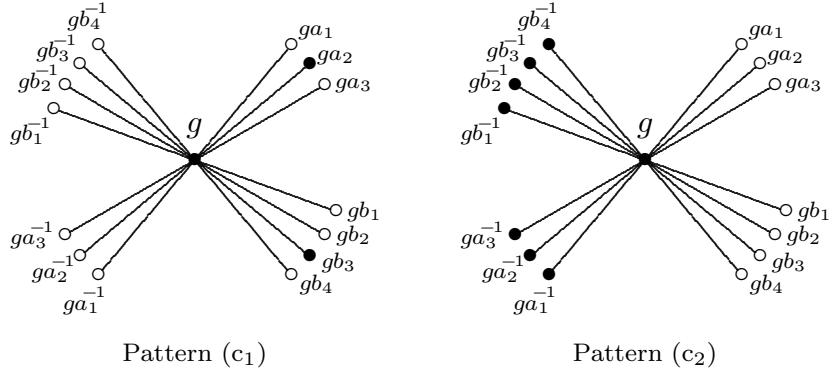
\newcount \ax \newcount \bx

\def\bull#1 #2 #3 #4 <#5>{%
  \ax = #1 \multiply \ax by 96 \divide \ax by 100
  \bx = #2 \multiply \bx by 96 \divide \bx by 100
  \put {\ifnum #3 = 0 $\circ$\else$\bullet$\fi} at #1 #2
  \plot 0 0 {\number\ax} {\number\bx} /
  \put {$\scriptstyle g\kern-0.5pt#4$} <#5> at #1 #2}

\def\iinv{\kern-6pt{\scriptscriptstyle^{^{-\kern-1pt1}}}}

\bigskip \noindent \hfill
\beginpicture
\setcoordinatesystem units <0.0020truecm, 0.0020truecm> point at 3000 0
\setplotarea x from -1000 to 1000, y from -1000 to 1000

\put {$\bullet$} at 0 0
\put {$g$} <0pt,10pt> at 0 0

\bull    866     500    0 {a_3} <9pt,0pt>
\bull    766     642    1 {a_2} <8pt,2pt>
\bull    642     766    0 {a_1} <8pt,4pt>

\bull    939    -342    0 {b_1} <9pt,0pt>
\bull    866    -500    0 {b_2} <9pt,0pt>
\bull    766    -642    1 {b_3} <9pt,-2pt>
\bull    642    -766    0 {b_4} <8pt,-3pt>

\bull   -866    -500    0 {a_3\iinv} <-10pt,0pt>
\bull   -766    -642    0 {a_2\iinv} <-8pt,-4pt>
\bull   -642    -766    0 {a_1\iinv} <-4pt,-8pt>

\bull   -939     342    0 {b_1\iinv} <-7pt,-7pt>
\bull   -866     500    0 {b_2\iinv} <-8pt,2pt>
\bull   -766     642    0 {b_3\iinv} <-7pt,5pt>
\bull   -642     766    0 {b_4\iinv} <-3pt,9pt>

\put{{\eightpoint Pattern (c$_1$)}} at 0 -1300

\setcoordinatesystem units <0.0020truecm, 0.0020truecm> point at 000 000

\put {$\bullet$} at 0 0
\put {$g$} <0pt,10pt> at 0 0

\bull    866     500    0 {a_3} <9pt,0pt>
\bull    766     642    0 {a_2} <8pt,2pt>
\bull    642     766    0 {a_1} <8pt,4pt>

\bull    939    -342    0 {b_1} <9pt,0pt>
\bull    866    -500    0 {b_2} <9pt,0pt>
\bull    766    -642    0 {b_3} <9pt,-2pt>
\bull    642    -766    0 {b_4} <8pt,-3pt>

\bull   -866    -500    1 {a_3\iinv} <-10pt,0pt>
\bull   -766    -642    1 {a_2\iinv} <-8pt,-4pt>
\bull   -642    -766    1 {a_1\iinv} <-4pt,-8pt>

\bull   -939     342    1 {b_1\iinv} <-7pt,-7pt>
\bull   -866     500    1 {b_2\iinv} <-8pt,2pt>
\bull   -766     642    1 {b_3\iinv} <-7pt,5pt>
\bull   -642     766    1 {b_4\iinv} <-3pt,9pt>

\put{{\eightpoint Pattern (c$_2$)}} at 0 -1300

\endpicture
\hfill\null
\bigskip\bigskip

\itemitem{(c$_1$)} there exists a unique $i\leq n$ and a unique $j\leq m$,
such that $ga_i$ and $gb_j$ lie in $\xi$, and for every $i$ and $j$,  none of
$ga_i\inv$ or $gb_j\inv$ lie in $\xi$,

\itemitem{(c$_2$)} for every $i$ and $j$, none of $ga_i$ or $gb_j$ lie in $\xi$,
and for every $i$ and $j$, all of $ga_i\inv$ and $gb_j\inv$ lie in $\xi$.

\bigskip

Having completed the description of $\OR$, the partial action of
$\Fmn$ is now easy to describe: for each $g\in \Fmn$ we put
  $$
  \OR_g = \big\{\xi\in \OR : g\in\xi\big\},
  $$
  and we let
  $$
  \tu_g: \OR_{g\inv} \to \OR_g,
  $$
  be given by $\tu_g(\xi) = g\xi = \{gh: h\in\xi\}$.

\bigskip

In possession of the proper notation we may now also describe the isomorphism
$\Psi$ mentioned in
\lcite{\IsoOnmProd}.  It is characterized by the fact that
  $$
  \Psi(\s(g)) = 1_g\univ\delta_g
  \for g\in \Fmn,
  \eqmark DescribePsi
  $$
  where $1_g\univ$ refers to the characteristic function of the clopen set $\OR_g\subseteq\OR$.

In what follows we will concentrate ourselves in studying the above partial
action of $\Fmn$ as well as the structure of $\Omn$ based on its crossed product
description.

\section Dynamical systems of type $(m,n)$

In this section we will study pairs of compact Hausdorff topological spaces
$(X,Y)$ such that
  $$
  X = \bigcup_{i=1}^n H_i = \bigcup_{j=1}^m V_j,
  $$
  where the $H_i$ are pairwise disjoint clopen subsets of $X$, each of which is
homeomorphic to $Y$ via given homeomorphisms
  $
  h_i: Y \to H_i.
  $
  Likewise we will assume that the $V_i$ are pairwise disjoint clopen subsets of $X$, each of which is
homeomorphic to $Y$ via given homeomorphisms
  $
  v_i: Y \to V_i.
  $
  See diagram \lcite{\DiagramOne}.

\definition We will refer to the quadruple
$\big(X,Y,\{h_i\}_{i=1}^n,\{v_j\}_{j=1}^m\big)$ as an \"{\mnsys}.

As an example, consider the situation  in which $Y\univ$ is the subset of
$\OR$ consisting of all the $\xi$ relative to which the configuration at $g=1$ follows
pattern ($c_2$).  Equivalently
  $$
  Y\univ = \big\{\xi \in \OR: a_i\inv,b_j\inv\in\xi, \hbox{ for all $i$ and $j$}\big\}.
  $$
  Let $X\univ$ be the complement of $Y\univ$ relative to $\OR$, and put
  $$
  h_i\univ:\xi\in Y\univ \mapsto a_i\xi\in X\univ
  \and
  v_j\univ:\xi\in Y\univ \mapsto b_j\xi \in X\univ.
  $$
  We leave it for the reader to verify that this provides an example of an
{\mnsys}.

\definition \label StdDynSys
  The system
  $\big(X\univ,Y\univ,\{h_i\univ\}_{i=1}^n,\{v_j\univ\}_{j=1}^m\big)$
  described above
  will be referred to as the \"{standard $(m,n)$--dynam\-ical system}.

It is our next immediate goal to prove that the standard {\mnsys} possesses a
universal property.  We thus fix, throughout, an arbitrary {\mnsys}
  $$
  \big(X,Y,\{h_i\}_{i=1}^n,\{v_j\}_{j=1}^m\big).
  $$

Our  goal will be to prove that there exists a unique map
  $$
  \gamma : X \mathop{\buildrel \cdot \over \cup}Y \to \OR
  $$
  such that
$\gamma(Y)\i Y\univ$,
$\gamma(X)\i X\univ$,
$\gamma\compos h_i = h_i\univ\compos\gamma$, and
$\gamma\compos v_j = v_j\univ\compos\gamma$.

We shall initially construct a partial action of $\Fmn$ on the topological disjoint
union
  $$
  \Omega := X \mathop{\buildrel \cdot \over \cup}Y.
  $$
  For this consider the inverse semigroup $\IO$ formed by all homeomorphisms
between clopen subsets of $\Omega$.  Evidently the $h_i$ and the $v_j$ are
elements of $\IO$.  Next consider the unique map
  $$
  \theta: \Fmn \to \IO
  $$
  such that
  $$
  \theta(a_i^{\pm 1}) = h_i^{\pm 1}, \quad
  \theta(b_j^{\pm 1}) = v_j^{\pm 1},
  $$
  and such that for each $g\in\Fmn$, written in reduced
  form\fn{That is, each $x_k$ is either $a_i^{\pm 1}$ or $b_j^{\pm 1}$, and
$x_{k+1}\neq x_k\inv$.}
  $$
  g = x_1x_2\ldots x_p,
  $$
  one has that
  $$
  \theta(g) = \theta(x_1)\theta(x_2)\ldots\theta(x_p).
  $$

  \state Proposition
  $\theta$ is a partial action of\/ $\Fmn$ on $\Omega$.

  \Proof
  It is not hard to prove this fact from scratch.  Alternatively one may
deduce it from known results as follows:
  using \cite[1.1]{\Pat},  one may faithfully represent $\IO$ as an inverse
semigroup of partial isometries on a Hilbert space.  Applying \cite[5.4]{\ortho}
we then conclude that there exists a unique semi-saturated partial representation
of $\Fmn$ in $\IO$, sending the $a_i$ to $h_i$, and the $b_j$ to $v_j$.  Evidently
this partial representation coincides with $\theta$,  and hence we conclude that
$\theta$ is a partial representation. Therefore,  for
every $g,h\in \Fmn$ one has that
  $$
  \theta_g \theta_h =
  \theta_g\theta_h\theta_{h\inv}\theta_h =
  \theta_{gh}\theta_{h\inv}\theta_h =
  \theta_{gh}\theta_h\inv\theta_h \subseteq
  \theta_{gh},
  $$
  meaning that $\theta_{gh}$ is an extension of $\theta_g \theta_h$, a property
that characterizes partial actions.
  \endProof

We may then form the crossed product $\CP$.  Given $g\in
\Fmn$, denote by $\Omega_g$ the range of $\theta_g$.  Since $\theta_g$ lies in
$\IO$, we have that its range is clopen.  So the characteristic function of
$\Omega_g$, which we shall denote by $1_g$, is a continuous function on
$\Omega$.

\state Proposition \label GeneralSemisatRep
  The map
  $$
  \rho:g\in \Fmn \mapsto  1_g\delta_g \in \CP,
  $$
  is a semi-saturated partial representation,  and moreover the elements
  $$
  s'_i:= \rho(a_i) \and t'_j:= \rho(b_j)
  $$
  satisfy relations $(\R)$.

  \Proof Given $g, h\in \Fmn$ we have
  $$
  \rho(g)\rho(h) =
  (1_g\d_g)(1_h\d_h) =
  \theta_g\big(\theta_{g\inv}(1_g)1_h\big)\d_{gh} =
  \theta_g\big(1_{g\inv}1_h\big)\d_{gh} =
  1_g1_{gh}\d_{gh}.
  \subeqmark RhoGH
  $$
  Therefore
  $$
  \rho(g)\rho(h)\rho(h\inv) =
  (1_g1_{gh}\d_{gh})(1_{h\inv}\d_{h\inv}) =
  \theta_{gh}\big(\theta_{(gh)\inv} (1_g1_{gh})1_{h\inv}\big) \d_g \$=
  \theta_{gh}\big(1_{h\inv}1_{(gh)\inv}1_{h\inv}\big) \d_g =
  \theta_{gh}\big(1_{(gh)\inv}1_{h\inv}\big) \d_g.
  $$
  On the other hand
  $$
  \rho(gh)\rho(h\inv) =
  (1_{gh}\d_{gh})(1_{h\inv}\d_{h\inv}) =
  \theta_{gh} \big(  \theta_{(gh)\inv}(1_{gh})1_{h\inv}\big)\d_g \$=
  \theta_{gh}\big(1_{(gh)\inv}1_{h\inv}\big) \d_g,
  $$
  which coincides with the above and hence proves that
  $
  \rho(g)\rho(h)\rho(h\inv) =
  \rho(gh)\rho(h\inv).
  $
  We leave it for the reader to prove that $\rho(g\inv) = \rho(g)^*$, after
which the verification that $\rho$ is a partial representation will be concluded.

Addressing semi-saturatedness, let $g,h\in \Fmn$ be such that $|gh|=|g|+|h|$.
This means that the reduced form of $gh$ is precisely the concatenation of the
reduced forms of $g$ and $h$, and hence we see that $\theta_{gh} =
\theta_g\compos\theta_h$.  In particular this implies that these two partial
homeomorphisms have the same range.  Therefore
  $$
  ran(\theta_g\compos\theta_h) = \theta_g(\Omega_{g\inv}\cap \Omega_h) =
  \Omega_g \cap \Omega_{gh}
  $$
  coincides with the range of $\theta_{gh}$,  which is $\Omega_{gh}$.  Having
concluded that $\Omega_g \cap \Omega_{gh}=\Omega_{gh}$, we deduce that
  $$
  1_g1_{gh} = 1_{gh}.
  $$
  Employing \lcite{\RhoGH} we then deduce that
  $$
  \rho(g)\rho(h) = 1_{gh}\d_{gh} = \rho(gh),
  $$
  proving that $\rho$ is semi-saturated.

Finally we leave it for the reader to prove that ${s'_i}^*s'_i$ and ${t'_j}^*t'_j$ coincide with
the characteristic function of $Y$, that $s'_i{s'_i}^*$ is the characteristic
function of $H_i$ (the range of $h_i$) and that $t'_j{t'_j}^*$ is the characteristic
function of $V_j$ (the range of $v_j$).  The checking of relations $(\R)$ now
becomes straightforward.
  \endProof

We may of course apply the above result for the standard {\mnsys}
  (see \lcite{\StdDynSys}),
  and hence there is a semi-saturated partial representation
  $$
  \rho\univ:\Fmn \to  \MCP,
  \eqmark DefineRhou
  $$
  given by
  $
  \rho\univ(g) = 1_g\univ\delta_g,
  $
  for every $g$ in $\Fmn$.
  With this notation \lcite{\DescribePsi} simply says that $\Psi\circ\s = \rho\univ$.

As another consequence of \lcite{\GeneralSemisatRep} and
\lcite{\FirstUnivResult} we have that there exists a *-homo\-morphism
  $$
  \Phi: \Omn \to \CP,
  $$
  such that $\rho = \Phi\circ\s$.

Wrapping up our previous results we obtain the commutative diagram:

\noindent \hfill
\beginpicture
\setcoordinatesystem units <0.0020truecm, 0.0020truecm> point at 3000 0
\setplotarea x from -700 to 2500, y from -1200 to 1300

\put{$C(\OR)\rtimes\Fmn$} at 1000 800
\arrow  <0.15cm> [0.25,0.75] from 1000 200 to 1000 600 \put{$\Psi$} at 1150 375
\put{$\Omn$} at 1000 000
\arrow  <0.15cm> [0.25,0.75] from 1000 -200 to 1000 -600 \put{$\Phi$} at 1150 -375
\put{$C(\Omega)\rtimes_\theta\Fmn$} at 1000 -800

\put{$\Fmn$} at -800 000

\arrow  <0.15cm> [0.25,0.75] from -650 150 to 180 650 \put{$\rho\univ$} at -200 580

\arrow  <0.15cm> [0.25,0.75] from -400 0 to 600 0 \put{$\s$} at 100 100

\arrow  <0.15cm> [0.25,0.75] from -650 -150 to 180 -650 \put{$\rho$} at -200 -570

\endpicture
\hfill\null

Observing that the correspondence
  $$
  f\in C(\Omega) \mapsto f\delta_1\in \CP
  $$
  is an embedding,   we will henceforth identify  $C(\Omega)$ with its image
within $\CP$ without further notice,  and similarly for $C(\OR)$.

\state Proposition  If\/ $\Gamma$ is defined as the composition $\Gamma := \Phi\Psi\inv$,
then $\Gamma \big(C(\OR)\big)\subseteq C(\Omega)$.

\Proof
  Since we will be dealing with two different dynamical systems here we will
insist in the convention (already used above) that
  $1_g$ denotes the characteristic function of $\Omega_g$, reserving
  $1\univ_g$ for the characteristic function of $\OR_g$.
  For each $g\in\Fmn$ we have that
  $$
  \rho\univ(g)\rho\univ(g\inv)=
  (1\univ_g\d_g) (1\univ_{g\inv}\d_{g\inv}) =
  1\univ_g\d_1 = 1\univ_g,
  $$
  and similarly
  $
  \rho(g)\rho(g\inv) = 1_g.
  $
  Since $\Gamma \circ \rho\univ = \rho$, we deduce that
  $$
  \Gamma(1\univ_g) = 1_g.
  \subeqmark GamaOneToOne
  $$

  It is easy to see that the set $\{1_g\univ:g\in\Fmn\}$ separates points of $\OR$, and
hence by the Stone-Weierstrass Theorem, the closed *-subalgebra it generates
coincides with $C(\OR)$.  So  the result follows from \lcite{\GamaOneToOne}.
  \endProof

As a consequence of the last result we see that there exists a unique continuous map
  $$
  \gamma : \Omega \to \OR,
  \eqmark UniqueGamma
  $$
  such that $\Gamma(f) = f\circ\gamma$, for every $f\in C(\OR)$.

\state Theorem \label UniversalPropertyOfSystem
  The standard {\mnsys} is universal in the following sen\-se: given any {\mnsys}
  $$
  \big(X,Y,\{h_i\}_{i=1}^n,\{v_j\}_{j=1}^m\big),
  $$
  there exists a unique continuous map
  $$
  \gamma : \Omega = X \mathop{\buildrel \cdot \over \cup}Y \to \OR,
  $$
  such that
  \izitem
  \zitem $\gamma(Y)\i Y\univ$,
  \zitem $\gamma(X)\i X\univ$,
  \zitem $\gamma\compos h_i = h_i\univ\compos\gamma$,
  \zitem $\gamma\compos v_j = v_j\univ\compos\gamma$.

\Proof  Regarding existence we will prove that the map $\gamma$ constructed in
\lcite{\UniqueGamma} satisfies the above properties.
Notice that
  $
  1_{a_1\inv}
  $
  is the characteristic function of the domain of $\theta(a_1)$ $(=h_1)$, namely
$Y$.  Similarly $1\univ_{a_1\inv}$ is the characteristic function of $Y\univ$.
Applying \lcite{\GamaOneToOne} to $g=a_1\inv$ we get
  $
  \Gamma(1_{Y\univ}) = 1_Y,
  $
  or equivalently
  $$
  1_{Y\univ}\circ\gamma = 1_Y.
  $$
  For $x\in\Omega$ this says that   $x\in Y$ iff $\gamma(x)\in Y\univ$, thus
proving  both (i) and (ii).

  Given $g\in \Fmn$, and $f\in C_0(\Omega\univ_{g})$, one may prove by direct
computation that
  $$
  \rho\univ(g\inv) f \rho\univ(g) =
  f\circ \theta\univ_{g},
  $$
  and similarly for $f\in C_0(\Omega_{g})$.
  So
  $$
  f\circ \theta\univ_{g} \compos \gamma =
  \Gamma\big(  f\circ \theta\univ_{g} \big) =
  \Gamma\big(  \rho\univ(g\inv) f \rho\univ(g)  \big) \$=
  \rho(g\inv) \Gamma(f) \rho(g) =
  \Gamma(f) \compos \theta_{g} =
  f\compos \gamma \compos \theta_{g}.
  $$
  Since $f$ is arbitrary it follows that
  $
  \theta\univ_{g} \compos \gamma = \gamma \compos \theta_{g}.
  $
  Point (iii) then follows by plugging $g=a_i$, while (iv) follows with $g=b_j$.

Addressing the uniqueness of $\gamma$, suppose one is given another map
  $$
  \gamma':\Omega \to\OR
  $$
  satisfying (i-iv).  Then it is clear that $\gamma'$  is covariant for the
corresponding partial actions of $\Fmn$ on $\Omega$ and $\OR$.
Letting
  $$
  \pi:f\in C(\OR) \mapsto f\compos \gamma'\in C(\Omega),
  $$
  one may easily prove that the pair $(\pi, \rho)$ is a covariant representation of the
partial dynamical system $(\OR,\tu,\Fmn)$ in $\CP$.  Using \cite[1.3]{\topfree}
we conclude that there exits a *-homomorphism
  $$
  \pi\times\rho: \MCP \to \CP,
  $$
  such that $(\pi\times\rho)(f) = f\circ\gamma'$, for every $f\in C(\OR)$, and such that
  $$
  (\pi\times\rho)\circ\rho\univ = \rho.
  \subeqmark GammaPrimeOnRho
  $$

Since the range of $\s$ generates $\Omn$, and since $\Psi$ is an isomorphism, we
deduce that the range of $\rho\univ$ generates $\MCP$.  We then conclude from
\lcite{\GammaPrimeOnRho} that
  $$
  (\pi\times\rho)\circ\rho\univ = \Gamma\circ\rho\univ,
  $$
  and hence that $\pi\times\rho=\Gamma$,  which in turn implies that $\gamma'=\gamma$.
\endProof

We shall next discuss the existence of fixed points in the universal
{\mnsys}.

\state Proposition \label FixedPoint
  If $n\geq m\geq2$, then there is a point $y$ in $Y\univ$ such that
  $$
  (v\univ_1)\inv h\univ_1(y) = y = (v\univ_2)\inv h\univ_2(y).
  $$

\Proof In order to prove the statement it is enough to show that there exists
some {\mnsys}
  $
  \big(X,Y,\{h_i\}_{i=1}^n,\{v_j\}_{j=1}^m\big),
  $
  and a point $y\in Y$ such that
  $$
  v_1\inv h_1(y) = y = v_2\inv h_2(y).
  $$
  By \lcite{\UniversalPropertyOfSystem},  the image of $y$ in $Y\univ$ under $\gamma$ will clearly
satisfy the required conditions.

We shall introduce another convenient variable by putting
  $$
  p := n-m+1.
  $$
  Let $Y=\{1,2,\ldots,p\}^\N$, with the product topology, and let $X$ be given
as the disjoint union of $m$ copies of $Y$.  To be precise,
  $$
  X = \{1,2,\ldots,m\}\times Y.
  $$
  For every $i=1, \ldots,  m$,  we define
  $$
  h_i:y\in Y\mapsto (i, y)\in X,
  $$
  and let us now define the $v_j$ via a process that is not as symmetric as
above.  For $j\leq m-1$, we put
  $$
  v_j:y\in Y\mapsto (j, y)\in X,
  \subeqmark AssymDefOne
  $$
  so that $v_j=h_j$, for all $j$'s considered so far.  In order to define the
remaining $v_j$'s, namely for $j$ of the form
  $$
  j = m-1 + k, \hbox{ with } k=1,\ldots,p,
  $$
  we let
  $$
  v_{m-1 + k}(y) = (m, ky) \for y\in Y,
  \subeqmark AssymDefTwo
  $$
  where ``$ky$" refers to the infinite sequence in $\{1,2,\ldots,p\}^\N$
obtained by preceding $k$ to $y$.
  The easy task of checking that the above does indeed gives an {\mnsys} is left for the reader.

We claim that the point $y =(1, 1, 1, 1,\ldots)$ satisfies the required conditions.
On the one hand we have the elementary calculation
  $$
  v_1\inv h_1(y) = v_1\inv(1,y) = y,
  $$
  where we are using the hypothesis that $m\geq2$, to guarantee that the
definition of $v_1$ is given by \lcite{\AssymDefOne} rather than by \lcite{\AssymDefTwo}.

If $m\geq3$, the same easy computation above yields
  $
  v_2\inv h_2(y) = y,
  $
  and the proof would be complete, so let us assume that $m=2$.
  Under this condition notice that $2=m-1+k$, with $k=1$, so
  $v_2$ is defined by \lcite{\AssymDefTwo}, and hence
  $$
  v_2(y) = (2,ky) = \big(2, k(1,1,1\ldots)\big)=
  \big(2, (1,1,1\ldots)\big) = h_2(y),
  $$
  whence $v_2\inv h_2(y) = y$,  and the claim is proven.

  As already mentioned, $\gamma(y)$ is then the element of $Y\univ$ satisfying the
requirements.
  \endProof

  % **********************************************************************************
  % The following section, typed by Pere, is inside a group since the conversion
  % from LaTeX to TeX required special macros and fonts.

  \begingroup

  \def\text#1{\hbox{#1}}
  \def\colon{\mathrel{:}}

  \font\tenmfrk   = eufm10
  \font\eightmfrk = eufm8
  \font\sixmfrk   = eufm6
  \def \tt {\fam \ttfam \eightmfrk}\relax
  \textfont \ttfam = \tenmfrk
  \scriptfont \ttfam = \eightmfrk
  \scriptscriptfont \ttfam = \sixmfrk

  \def \mathfrak#1{{\tt #1}}
  \def \mathbb#1{{\bf #1}}
  \def \mathcal#1{{\cal #1}}
  \def \\{\hfill\pilar{12pt}\cr}
  \def\eqalign#1{ \vcenter{ \ialign{ \pilar{10pt} $##$&\ $##$\hfil\crcr #1\crcr}}}
  \def\refer #1:#2;{%
    \ifundef {#2} UNDEFINED LABEL (#2), EXITING.  \end \fi
    \csname #2\endcsname \relax}
  \def\ref#1{\refer #1;}
  \def\parref#1{(\refer #1;)}
  \def\subset{\subseteq}

\section Configurations and functions

  The purpose of this section is to give a description of the space $Y\univ $ of
configurations of pattern $(c_2)$ at $1$. This will be done in terms of certain
functions, which we are now going to describe.

Set $Z_0:=\{ a_1,\dots ,a_n\}$, $Z_1:=\{ b_1,\dots , b_m\}$ and $E=
Z_0\sqcup Z_1$. We will denote the elements of $E^r$ as words
$\alpha =e_1e_2\cdots e_r$ in the alphabet $E$. Set $E^+:= \bigsqcup
_{r=1}^{\infty} E^r$.   For $\alpha =e_1e_2\cdots e_r\in E^+$ define
the color of $\alpha$ as $c(\alpha )= 1-i$, if $e_r\in Z_i$. We
consider the compact Hausdorff space
$$Z:=\prod _{\alpha \in E^+} Z_{c(\alpha)},$$
where each $Z_i$ is given the discrete topology and $Z$ is endowed
with the product topology. Elements of $Z$ will be interpreted as
functions $f\colon E^+\to E$ such that $f(\alpha)\in Z_{c(\alpha)}$
for all  $\alpha \in E^+$.

Let $D$ be the subspace of $Z$ consisting of the functions $f$ such
that the following properties (*) and (**) hold for all $\alpha \in
E^+\sqcup \{\cdot \}$, all $e\in E$ and all $\beta \in E^+$:

\medskip

$\hskip0.05cm $ (*) $\, f(\alpha ef(\alpha e))=e$.

(**) $\, f(\alpha ef(\alpha e)\beta )= f(\alpha \beta)$.

\medskip

Observe that, for $e\in E$ and $\alpha \in E^+\sqcup \{\cdot \}$, we
have $\,\,  e \in Z_i \iff f(\alpha e)\in Z_{1-i}$. It is easy to
show that $D$ is a closed subspace of $Z$, and thus $D$ is a compact
Hausdorff space with the induced topology.

Our aim in this section  is to show the following result:

\state Theorem \label FuncEqualConfig
  There is a canonical homeomorphism $D\cong Y\univ$.

To show this we need some preliminaries.

\definition
  A {\it partial $E$-function} is a family
  $(\Omega_1,f_1),(\Omega_2,f_2),\dots , (\Omega_r,f_r)$, for some $r\ge 1$,
satisfying the following relations:

\initem
  \nitem $\Omega _1 = E$, and $f_1\colon E\to E$ is a function such that
$f_1(e)\in Z_{c(e)}$ for all $e\in E$.
  \nitem For each $i=1, \dots ,r$, $$\Omega _i=\{ x_1x_2\cdots x_i\in E^i \mid
x_{j+1}\ne f_j(x_1x_2\cdots x_j) \text{ for } j=1,\dots , i-1 \}, $$ and
$f_i\colon \Omega _i \to E$ is a function such that $f_i(\alpha)\in
Z_{c(\alpha)}$ for all $\alpha \in \Omega _i$.
  \medskip\noindent
An {\it $E$-function} is an infinite sequence
$(\Omega_1,f_1),(\Omega_2,f_2),\dots , $ satisfying the above
conditions for all indices.

It is quite clear that any partial
$E$-function can be extended (in many ways) to an $E$-function.

\state Lemma \label Efunctions
  Given an $E$-function $(\Omega_1,f_1),(\Omega _2,f_2),\dots$, there is a
unique function $f\in D$ such that $f(\alpha)=f_i(\alpha)$ for $\alpha \in
\Omega _i$ and all $i\in \N$. Therefore $D$ can be identified with the space of
all $E$-functions. Moreover a basis for the topology of $D$ is provided by the
partial $E$-functions by the rule: $$\mathfrak f=((\Omega_1, f_1), (\Omega
_2,f_2), \dots , (\Omega _r,f_r))\longmapsto U_{\mathfrak f}\, ,$$ where
$U_{\mathfrak f}=\{f\in D\mid f \,\, {\rm extends } \, \,\mathfrak f \}$.

  \Proof
Let $(\Omega_1,f_1),(\Omega _2,f_2),\dots , $ be an $E$-function. We
have to construct an extension of it to a function  $f\colon E^+\to
E$ such that $f\in D$. It will be clear from the construction that
$f$ is unique.

Note that $f(ef(e))$ must be equal to $e$ for $e\in E$ by condition
(*). This, together with the extension property determines
completely $f$ on $E^{\le 2}$. Assume that $f$ has been defined on
$E^{r-1}$ for some $r\ge 3$. Then we define $f$ on $E^r$ as follows:
First $f(\alpha )= f_r (\alpha)$ if $\alpha\in \Omega_r$. If $\alpha
=x_1x_2\cdots x_r\notin \Omega _r$, there are various possibilities,
that we are going to consider:

If $x_2=f(x_1)$, then we set
$$f(x_1f(x_1)x_3\cdots x_r)= f(x_3\cdots x_r).$$
Observe that this is forced by condition (**).

Analogously, if $x_{j+1}\ne f(x_1\cdots x_j)$ for $j=1,\dots ,i-1$
and $x_{i+1}=f(x_1x_2\cdots x_i)$ for some $i<r-1$, define
$$f(x_1x_2\cdots x_if(x_1\cdots x_i)x_{i+2}\cdots x_r)=f(x_1x_2\cdots x_{i-1}x_{i+2}\cdots x_r).$$
Also we have here that this is forced by (**).

Finally if $x_{j+1}\ne f(x_1\cdots x_j)$ for $j=1,\dots ,r-2$ and
$x_{r}=f(x_1x_2\cdots x_{r-1})$, define
$$f(x_1x_2\cdots x_{r-1}f(x_1\cdots x_{r-1}))=x_{r-1}.$$
Note that this is forced by (*).

We obtain a map $f\colon E^+\to E$ such that $f(\alpha )\in
Z_{c(\alpha)}$ for all $\alpha\in E^+$. We have to check conditions
(*) and (**).

For (*), let $\alpha\in E\cup \{\cdot\}$ and $e\in E$. We will check
that $f(\alpha ef(\alpha e))= e$ by induction on $| \alpha |$. If
$\alpha = \cdot $ then we have that $f(ef(e))=e$ by construction.
Suppose that the equality holds for words of length $r$ and let
$\alpha $ a word of length $r+1$. Write $\alpha = x_1x_2 \cdots
x_{r+1}$. Assume that, for $1\le i\le r$, we have that $x_{j+1}\ne
f(x_1\cdots x_j)$ for all $j=1,\dots ,i-1$, and that
$x_{i+1}=f(x_1x_2\cdots x_i)$. Then we have

$$\eqalign{
f(\alpha ef(\alpha e)) & = f(x_1x_2\cdots x_if(x_1x_2\cdots
x_i)x_{i+2}\cdots x_{r+1}ef(\alpha e))\\
& = f(x_1\cdots x_{i-1}x_{i+2}\cdots x_{r+1}ef(\alpha e))\\
& = f(x_1\cdots x_{i-1}x_{i+2}\cdots x_{r+1}ef(x_1x_2\cdots
x_if(x_1x_2\cdots x_i)x_{i+2}\cdots x_{r+1}e))\\
& = f(x_1\cdots x_{i-1}x_{i+2}\cdots x_{r+1}ef(x_1x_2\cdots
x_{i-1}x_{i+2}\cdots x_{r+1}e))\\
& = e}$$

\medskip\noindent
where we have used the induction hypothesis in the last step.

Assume that $x_{j+1}\ne f(x_1\cdots x_j)$ for all $j=1,\dots ,r$,
and that $e=f(x_1\cdots x_rx_{r+1})$. Then we have

$$\eqalign{
f(\alpha ef(\alpha e)) & = f(x_1\cdots x_rx_{r+1}f(x_1\cdots
x_rx_{r+1})f(\alpha e))\\
& = f(x_1\cdots x_rf(\alpha e))\\
& = f(x_1\cdots x_{r}f(x_1x_2\cdots x_rx_{r+1}f(x_1x_2\cdots x_rx_{r+1})))\\
& = f(x_1\cdots x_{r}x_{r+1})\\
& = e.}$$

\medskip
Finally if $x_{j+1}\ne f(x_1\cdots x_j)$ for all $j=1,\cdots ,r$ and
$e\ne f(x_1x_2\cdots x_{r+1})$ then we have
$$f(x_1\cdots x_{r+1}ef(x_1\cdots x_{r+1}e))=e$$
by the definition of $f$.

The checking of (**) is similar. We prove that $f(\alpha ef(\alpha
e)\beta)=f(\alpha \beta)$ by induction on $|\alpha |$. If $\alpha
=\cdot $ then $f(ef(e)\beta)= f(\beta)$ by definition of $f$.
Suppose (**) holds when the length of $\alpha$ is $\le r$ and set
$\alpha =x_1x_2\cdots x_{r+1}$.  Assume that, for $1\le i\le r$, we
have that  $x_{j+1}\ne f(x_1\cdots x_j)$ for all $j=1,\dots ,i-1$,
and that $x_{i+1}=f(x_1x_2\cdots x_i)$. Then we have

$$\eqalign{
f(\alpha ef(\alpha e)\beta) & = f(x_1x_2\cdots x_if(x_1x_2\cdots
x_i)x_{i+2}\cdots x_{r+1}ef(\alpha e)\beta)\\
& = f(x_1\cdots x_{i-1}x_{i+2}\cdots x_{r+1}ef(x_1\cdots x_{i-1}x_{i+2}\cdots x_{r+1}e)\beta)\\
& = f(x_1\cdots x_{i-1}x_{i+2}\cdots x_{r+1}\beta)\\
& = f(x_1\cdots x_{r+1}\beta)\\
& = f(\alpha \beta)}$$

\medskip\noindent
where we have used the induction hypothesis for the third equality.

Assume that $x_{j+1}\ne f(x_1\cdots x_j)$ for all $j=1,\dots ,r$,
and that $e=f(x_1\cdots x_rx_{r+1})$. Then we have that $f(\alpha
e)= f(x_1\cdots x_rx_{r+1}e)=x_{r+1}$ by definition of $f$, and so

$$\eqalign{
f(\alpha ef(\alpha e)\beta) & = f(x_1\cdots x_rx_{r+1}f(x_1\cdots
x_rx_{r+1})f(\alpha e)\beta)\\
& = f(x_1\cdots x_rf(\alpha e)\beta)\\
& = f(x_1\cdots x_{r}x_{r+1}\beta )\\
& = f(\alpha \beta).}$$

\medskip
Finally if $x_{j+1}\ne f(x_1\cdots x_j)$ for all $j=1,\cdots r$ and
$e\ne f(x_1x_2\cdots x_{r+1})$ then we have
$$f(x_1\cdots x_{r+1}ef(x_1\cdots x_{r+1}e)\beta)=f(x_1\cdots x_{r+1}\beta)=f(\alpha\beta)$$
by the definition of $f$.

Given the description of $D$ as a subspace of $Z=\prod _{\alpha \in
E^+}Z_{c(\alpha)}$, it is clear that the family $\{U_{\mathfrak
f}\mid \mathfrak f \text{ is a partial } E-\text{function}\}$ is a
basis for the topology of $D$.
  \endProof

\noindent {\it Proof of Theorem \parref{thm:FuncEqualConfig}.} We will
define mutually inverse maps $\varphi \colon Y\univ\to D$ and $\psi
\colon D \to Y\univ$.

Let $\xi\subset \Fmn$ be a configuration of pattern
$(c_2)$ at $1$. Then we have that $x^{-1}\in \xi$ for all $x\in E$.
Now the configuration at $x^{-1}$ must be of pattern $(c_1)$, so
that, for each $x\in E$ there is a unique $f_1(x)\in Z_{c(e)}$ such
that $x^{-1}f(x)\in \xi$. This defines a partial $E$-function
$(E,f_1)$. For each $x_1\in E$,  the configuration at
$x_1^{-1}f_1(x_1)$ must be of pattern $(c_2)$, so all words of the
form $x_1^{-1}f(x_1)x_2^{-1}$, with $x_2\ne f_1(x_1)$ must be in
$\xi$. In the next step we look at the configuration at vertices of
the form $x_1^{-1}f(x_1)x_2^{-1}$, where $x_2\ne f_1(x_1)$. Here the
configuration must be of pattern $(c1)$, so there is a unique
$f_2(x_1x_2)\in Z_{c(x_2)}=Z_{c(x_1x_2)}$ such that
$x_1^{-1}f_1(x_1)x_2^{-1}f_2(x_1x_2)\in \xi$. This gives us a
partial $E$-function $(\Omega _1,f_1), (\Omega _2, f_2)$, where of
course $\Omega _2=\{x_1x_2\in E^2\mid x_2\ne f_1(x_1)\}$. Proceeding
in this way we obtain an $E$-function $\varphi (\xi )=((\Omega _1,
f_1), (\Omega _2,f_2), \dots )$.

To define $\psi$  we just need to revert the previous process. Given
an $E$-function $$\mathfrak f= ((\Omega _1,f_1),
(\Omega_2,f_2),\dots ),$$ the configuration $\psi (\mathfrak f)$
consists of $1$ together with the elements of $\Fmn$ of
the form
$$x_1^{-1}f_1(x_1)x_2^{-1}f_2(x_1x_2)\cdots x_r^{-1}$$
and
$$x_1^{-1}f_1(x_1)x_2^{-1}f_2(x_1x_2)\cdots x_r^{-1}f_r(x_1x_2\cdots x_r). $$
where $x_1x_2\cdots x_r\in \Omega _r$, for $r\ge 1$. It is clear
that $\psi (\mathfrak f)$ is a configuration of pattern $(c_2)$ at
$1$, and that $\varphi $ and $\psi $ are mutually inverse maps.

Since both $D$ and $Y\univ$ are compact Hausdorff spaces, in order
to show that $\varphi$ is a homeomorphism it is enough to prove that
$\psi $ is an open map. Since the family $\{U_{\mathfrak f}\mid
\mathfrak f \text{ is a partial } E-\text{function}\}$ is a basis
for the topology of $D$ by Lemma
  \parref{lem:Efunctions},
  it is enough
to show that $\psi (U_{\mathfrak f})$ is an open subset of $Y\univ$
for every partial $E$-function $\mathfrak f$. Thus let $\mathfrak f
= ((\Omega _1,f_1), (\Omega _2,f_2),\dots , (\Omega_r,f_r))$ be a
partial $E$-function, and consider the set
$$T:= \{ x_1^{-1} f_1(x_1)x_2^{-1}f_2(x_1x_2)\cdots x_r^{-1} f_r (x_1x_2\cdots x_r)\mid
x_1x_2\cdots x_r \in \Omega _r\}.$$ By using the convexity of the
elements of $Y\univ$ it is straightforward to show that
$$\psi (U_{\mathfrak f} )= \{ \xi \in Y\univ \mid g\in \xi \quad \forall g\in T \}.$$
Since $T$ is a finite subset of $\Fmn$, we conclude that
$\psi (U_{\mathfrak f})$ is an open subset of $Y\univ$. This
concludes the proof of Theorem \parref{thm:FuncEqualConfig}.

It will be useful to get a detailed description of the action $\tu$ of
$\Fmn$ on $Y\univ $ in terms of the picture of $Y\univ$
using  $E$-functions (Theorem \ref{thm:FuncEqualConfig}).

\state Lemma \label actionForFuncts
  Let $$g= z_r^{-1}x_rz_{r-1}^{-1} x_{r-1}\cdots z_1^{-1}x_1 $$ be a reduced
word in $\Fmn$, where $x_1,\dots ,x_r,z_1,\dots ,z_r\in E$. Then
$\text{Dom}(\tu_g)=\emptyset $ unless $z_i\in Z_{c(x_i)}$ for all $i=1,\dots
,r$. Assume that the latter condition holds. Then the domain of $g$ is precisely
the set of all $E$-functions $\mathfrak f=(f_1,f_2,\dots )$ such that
$f_i(x_1\cdots x_i)=z_i$ for all $i=1,\dots ,r$, and the range of $g$ is the set
of those $E$-functions $\mathfrak h=(h_1,h_2,\dots , )$ such that
$h_i(z_rz_{r-1}\cdots z_{r-i+1})= x_{r-i+1}$ for all $i=1,\dots ,r$. Moreover
for $\mathfrak f \in \text{Dom}(\tu_g)$ let $\mathfrak h= {^{g}\mathfrak f}$ denote
the image of $\mathfrak f$ under the action of $g$. Then $\mathfrak h=((\Omega
_1,h_1),(\Omega _2, h_2),\dots )$ with
  $$
  h_{r+t}(z_rz_{r-1}\cdots z_1y_1y_2\cdots y_t)=
f_t(y_1\cdots y_t) \quad \text{if}\quad  z_r\cdots z_1y_1\cdots
y_t\in \Omega _{r+t}.
  \subeqmark actionOne
  $$
  Moreover, for $i=2,3,\dots , r$ and $z_rz_{r-1}\cdots z_iy_1\cdots
y_t\in \Omega _{r-i+1+t}$ with $z_{i-1}\ne y_1$,
  $$
  h_{r-i+1+t}(z_rz_{r-1}\cdots  z_i y_1y_2\cdots
y_t)= f_{i-1+t}(x_1\cdots x_{i-1}y_1y_2\cdots y_t),
  \subeqmark actionTwo
  $$
and for $y_1\cdots y_t\in \Omega _t$ with $y_1\ne z_r$ we have
  $$
  h_t(y_1\cdots y_t)= f_{r+t}(x_1x_2\cdots
x_ry_1y_2\cdots y_t)\, .
  \subeqmark actionThree
  $$

  \Proof
Suppose that $\text{Dom}(\tu_g)\ne \emptyset$ and take $\xi \in
\text{Dom}(\tu_g)$, with corresponding $E$-function $\mathfrak f$. Since
$1\in \xi$ we get that $g\in \xi$ and by convexity we get that
$z_r^{-1}x_rz_{r-1}^{-1}x_{r-1}\cdots z_i^{-1} x_i\in \xi$ for all
$i=1,\dots ,r$. We thus obtain that $h_i(z_rz_{r-1}\cdots
z_{r-i+1})= x_{r-i+1}$ for all $i=1,\dots ,r$, where $\mathfrak h=
{^{g}\mathfrak f}$. In particular it follows that $x_{r-i+1}\in
Z_{c(z_{r-i+1})}$ for $i=1,\dots ,r$, which is equivalent to $z_i\in
Z_{c(x_i)}$ for $i=1,\dots, r$. Moreover since $gx_1^{-1}f_1(x_1)\in
\xi$, we get $z_r^{-1}x_r\cdots x_2z_1^{-1}f_1(x_1)\in \xi$. Since
$z_1,f(x_1)\in Z_{c(x_1)}$ we get that $z_1=f_1(x_1)$. Similarly we
get that $f_i(x-1\cdots x_i)=z_i$ for all $i=1,\dots ,r$.

Conversely, assume that $z_i\in Z_{c(x_i)}$ for all $i=1,\dots ,r$.
Then there are infinitely many $E$-functions $\mathfrak f =(f_1,f_2,
\dots )$ such that $f_i(x_1x_2\cdots x_i)= z_i$ for all $i=1,\dots
,r$. Let $\mathfrak f$ be one of these functions. Then it is easy to
verify that ${^{g}\mathfrak f} $ is the $E$-function $\mathfrak h
=(h_1,h_2,\dots )$ determined by $h_i(z_rz_{r-1}\cdots
z_{r-i+1})=x_{r-i+1}$ for $i=1,\dots ,r$ and by the rules
\parref{eq:actionOne}, \parref{eq:actionTwo} and \parref{eq:actionThree}.
  \endProof

Recall the following definition from
\cite{\topfree}.

\definition \label DefTopFree
  Let $\theta $ be a partial action of a group $G$ on a compact Hausdorff space
$X$. The partial action $\theta$ is {\it topologically free} if for every $t\in
G\setminus \{1\}$, the set $F_t:=\{ x\in U_{t^{-1}}\mid \theta _t(x)=x\}$ has
empty interior.

\state Proposition \label Onmtopfree
  For $m,n\ge 2$, the action of\/ $\Fmn$ on $\OR$ is topologically
free.

  \Proof
Let $g\in \Fmn\setminus \{1\}$. Assume first that
$$g=z_r^{-1}x_rz_{r-1}^{-1} x_{r-1}\cdots z_1^{-1}x_1$$
is a reduced word with $x_1,\dots ,x_r,z_1,\dots ,z_r\in E$.
Obviously we may suppose that the domain of $g$ is non-empty, so
that $z_i\in Z_{c(x_i)}$ for $i=1,\dots ,r$ by Lemma
\parref{lem:actionForFuncts}.

Note that the domain of $g$ is contained in $Y\univ$. By Theorem
\parref{thm:FuncEqualConfig} we only have to show that for any partial
$E$-function $\mathfrak f$ there is an extension $\mathfrak f '$ of
$\mathfrak f $ such that ${^{g}\mathfrak f '}\ne \mathfrak f '$.
Obviously we can assume that $\mathfrak f = ((\Omega _1,f_1),(\Omega
_2, f_2),\dots , (\Omega _s,f_s))$, with $s> r$ and that
$U_{\mathfrak f}\cap \text{Dom}(\tu_g)\ne \emptyset$. Set $t=s-r$ and
choose $y_1,\dots ,y_t$ in $E$ such that $z_rz_{r-1}\cdots
z_1y_1\cdots y_t\in \Omega _s$. Select $y_{t+1}\in E$ such that
$y_{t+1}\ne f_s(z_r\cdots z_1y_1\cdots y_t)$. Since $m, n\ge 2$,
there exits $u\in Z_{c(y_{t+1})}$ such that $u\ne f_{t+1}(y_1\cdots
y_ty_{t+1})$. Define $f_{s+1}\colon \Omega _{s+1}\to E$ in such a
way that $f_{s+1}(z_1\cdots z_ry_1\cdots y_ty_{t+1})=u$, and
arbitrarily on the other elements of $\Omega _{s+1}$ subject to the
condition that $f_{s+1}(w_1\cdots w_{s+1})\in Z_{c(w_{s+1})}$. Then
$((\Omega _1,f_1),\dots , (\Omega _s,f_s), (\Omega _{s+1},f_{s+1}))$
is a partial $E$-function extending $\mathfrak f$. Extend this
partial $E$-function to an $E$-function $\mathfrak f'$. If
${^{g}\mathfrak f'}=\mathfrak f '$ then equation \parref{eq:actionOne}
gives
$$f_{s+1}(z_r\cdots z_1y_1\cdots y_ty_{t+1})= f_{t+1}(y_1\cdots
y_ty_{t+1}),$$ which contradicts our choice of $f_{s+1}(z_r\cdots
z_1y_1\cdots y_ty_{t+1})$.

We conclude that $U_{\mathfrak f}$ has points which are not fixed
points for $g$.

Now assume that
$$g=x_rz_{r-1}^{-1}x_{r-1}\cdots z_1^{-1}x_1 z_0^{-1}$$
is a reduced word in $\Fmn$, with $x_1,\dots
,x_r,z_0,\dots ,z_{r-1}\in E$. Write
  $$
  g':=x_rz_{r-1}^{-1}x_{r-1}\cdots z_1^{-1}x_1.
  $$
  Assume that $g\cdot
\xi =\xi$ for all $\xi\in V$, where $V$ is an open subset of $X$.
Then $(z_0^{-1}g')\cdot \xi '=\xi '$ for all $\xi '\in z_0^{-1}V$.
By the first part of the proof we get $z_0^{-1}g'=1$ and thus
$g=g'z_0^{-1}=1$, as desired.
  \endProof

As an easy consequence we obtain:

\state Corollary If $\rho$ is a representation of $\Omnr$ whose restriction to
$C(\OR)$ is injective, then $\rho$ itself is injective.

\Proof
  Follows immediately from \lcite{\Onmtopfree} and  \cite[2.6]{\topfree}.
  \endProof

Recall the following definition.

  \definition
  A $C^*$-algebra satisfies property (SP) (for {\it
small projections}) in case every nonzero hereditary
$C^*$-subalgebra contains a nonzero projection. Equivalently, for
every nonzero positive element $a$ in $A$ there is $x\in A$ such
that $x^*ax$ is a nonzero projection.

\state Theorem \label SPproperty
  For $m,n\ge 2$, the $C^*$-algebra $\mathcal O^r_{m,n}$ satisfies property
${\rm (SP)}$. More precisely, given a nonzero positive element $c$ in $\mathcal
O^r_{m,n}$, there is an element $x\in \mathcal O^r_{m,n}$ such that $x^*cx$ is a
nonzero projection in $C(\OR)$. In particular every nonzero ideal of $\mathcal O
^r_{m,n}$ contains a nonzero projection of $C(\OR)$.

  \Proof
This is well-known for the Cuntz algebras $\mathcal O_n$ so we may
assume that $m,n\ge 2$.

Let $c$ be a nonzero positive element in $\mathcal O^r_{m,n}$. Since
the canonical conditional expectation $E_r$ is faithful, we may
assume that $\| E_r(c)\|=1$. By Proposition \parref{prop:Onmtopfree} the
partial action of $\Fmn$ on $\OR$ is
topologically free. Hence, it follows from \cite[Proposition
2.4]{\topfree} that, given $1/4 >\epsilon >0$, there is an element $h\in
C(\OR)$ with $0\le h\le 1$ such that
  \initem
\nitem $\| hE_r(c)h\| \ge \| E_r(c)\| -\epsilon$,
\nitem $\|hE_r(c)h -hch\| \le \epsilon $.
  \medskip\noindent
  By \cite[Lemma 2.2]{\KR} there is a contraction $d$ in $\mathcal
O^r_{m,n}$ such that $d^*(hch)d= (hE_r(c)h-\epsilon )_+$, and so it
follows that $(hd)^*c(hd)$ is a nonzero positive element in
$C(\OR)$. Since $C(\OR)$ is an AF-algebra it has
property (SP) so there is an element $y$ in $C(\OR )$ such
that $y^*(hd)^*c(hd)y$ is a nonzero projection in $C(\OR)$.
Taking $x= hdy$, we get the result.
  \endProof

\endgroup

%=====================================================================
%=====================================================================
%=====================================================================

\section Exactness of the reduced cross-sectional C*-algebra of a Fell bundle

Recall from \cite[Section 2]{\infinoa} that the full (resp.~reduced)
crossed product may be defined as the full (resp.~reduced) cross
sectional C*-algebra of the semidirect product Fell bundle
\cite[2.8]{\tpa}.
  For this reason we shall now pause to prove some key results on Fell bundles
in support our study of $\Omn$.

We begin by discussing the notion of (minimal) tensor product of a C*-algebra by
a Fell bundle.  We refer the reader to \cite{\FD} for an extensive study of the
theory of Fell bundles.  We thank N.~Brown for an interesting conversation from
which some of the ideas pertaining to this tensor product arose.

\state Proposition \label TensorProdFBun
  Let $A$ be a C*-algebra and let $\B = \{B_g\}_{g\in G}$ be a Fell bundle over
a discrete group $G$.  Then
  \izitem
  \zitem
  There exists a unique collection of seminorms
$\{\|\cdot\|_g\}_{g\in G}$ on
the algebraic tensor products $A\.B_g$, such that
  $\|\cdot\|_1$ is the spacial (minimal) C*-norm on $A\. B_1$, and
the completions
  $$
  A\*B_g := \overline{A\. B_g}^{\|\cdot\|_g}
  $$
  become the fibers of a Fell bundle
  $\{A\*B_g\}_{g\in G}$,
  in which the multiplication and involution operations extend the following:
  $$
  \matrix {
  (a_1\* b_1,  a_2\* b_2) \in \big(A\. B_{g_1}\big) \times \big(A\. B_{g_2}\big)
  & \longmapsto &
  a_1 a_2 \*b_1 b_2 \in A\. B_{g_1g_2} \cr
  \hfill  \pilar{20pt}
  (a\* b) \in A\. B_g
  & \longmapsto &
  a^* \*b^* \in A\. B_{g\inv}. \hfill}
  $$
  \zitem  Denoting the resulting Fell bundle by $A\*\B$, there exists a (necessarily unique) *-isomor\-phism
  $$
  \phi:A\*\Cr\B \to \Cr{A\*\B},
  $$
  such that $\phi(a\*b_g) = a\*b_g$, whenever $a\in A$, and $b\in B_g$, for any
$g$ (the last two tensor product signs should be given the
appropriate and obvious meaning in each case).

\Proof
  In order to prove uniqueness, suppose that a collection of norms is given as
above.  Then,
for every $g\in G$,  and any $c\in A\.B_g$, one has that
  $$
  \|c\|_g^2 = \|c^*c\|_1.
  $$
  Since $c^*c\in A\.B_1$, and since the norm on $A\.B_1$ is assumed to be the
spacial norm,  uniqueness immediately follows.  As for existence, let
  $$
  \pi:A\to B(H)
  $$
  be a faithful representation of $A$ on a Hilbert space $H$, and let
  $$
  \rho: \bigcup_{g\in G} B_g \to B(K)
  $$
  be a representation (in the sense of \cite[2.2]{\amena}) of $\B$ on a Hilbert
space $K$, which is isometric on each $B_g$.  Such a representation may be
easily obtained by composing the natural inclusion maps $B_g\to \Cr\B$, which
are isometric by \cite[2.5]{\amena}, with any faithful representation of
$\Cr\B$.

  Consider the representations
  $$
  \matrix{
  \pi' = \pi\*1& : & A & \to & B(H\*K) \cr
  \pilar{20pt}
  \rho' = 1\*\rho& : & \bigcup_{g\in G} B_g & \to & B(H\*K),
  }
  $$
  and let
  $$
  C_g = \clspan\big( \pi'(A)\rho'(B_g) \big).
  $$
  It is then  easy to see that $C_gC_h\subseteq C_{gh}$, and $C_g^*\subseteq
C_{g\inv}$, for every $g,h\in G$.   So we may think of $\C = \{C_g\}_{g\in G}$ as
a Fell bundle over $G$,  with operations borrowed from $B(H\*K)$.

For each $g$ in $G$ consider the seminorm $\|\cdot\|_g$ on $A\.B_g$ obtained as
the result of composing the maps
  $$
  A\.B_g
  \buildrel \pi'\*\rho' \over \longrightarrow
  C_g
  \buildrel \|\cdot\| \over \longrightarrow
  {\bf R}.
  $$
  Evidently the completion of $A\.B_g$ under this seminorm is isometrically
isomorphic to $C_g$.  By \cite[3.3.1]{\Nate} we have that $\|\cdot\|_1$ is the
spatial norm and the remaining conditions in (i) may now easily be verified.

In order to prove (ii) we consider two other representations of our objects,
namely
  $$
  \matrix{
  \pi'' = \pi'\*1 = \pi\*1\*1& : & A & \to & \BigB \cr
  \pilar{20pt}
  \rho'' = \rho'\*\l = 1\*\rho\*\l& : & \bigcup_{g\in G} B_g & \to & \BigB,
  }
  $$
  where $\l$ is the regular representation of $G$, and for any given $b_g$ in
$B_g$, we put
  $$
  \rho''(b_g) = 1\*\rho(b_g)\*\l_g.
  $$

Observing that $\rho''$ is also isometric on each $B_g$, we see that
the closed
*-subalgebra of $\BigB$ generated by the range of $\rho''$ is isomorphic to
$\Cr\B$ by \cite[3.7]{\amena} (the faithful conditional expectation is just the
restriction to the diagonal).  Alternatively one may also deduce this from
\cite[3.4]{\exact}.

By \cite[3.3.1]{\Nate} one then has that $A\*\Cr\B$ is isomorphic to
the subalgebra of operators generated by $\pi''(A)\rho''(\B)$.
  For further reference let us observe that the present model of
$A\*\Cr\B$ within $\BigB$ is therefore generated by the set
  $$
  \{a\*(b_g\*\l_g): a\in A,\  g\in G,\ b_g\in B_g\}.
  $$

Observe that,  for each $g\in G$,  the map
  $$
  \sigma_g\ :\  x \in C_g \longmapsto x\*\l_g\in \BigB
  $$
  is an isometry and, collectively, they provide a representation of
  $\C$ in $\BigB$.

By the same reasoning employed above, based on \cite[3.7]{\amena} or \cite[3.4]{\exact}, we have that
$\Cr{\C}$ is isomorphic to the closed *-subalgebra of $\BigB$ generated by
the union of the ranges of all the $\sigma_g$.
  Therefore our model of  $\Cr{\C}$ within $\BigB$ is
generated by the set
  $$
  \{(a\*b_g)\*\l_g: a\in A,\  g\in G,\ b_g\in B_g\}.
  $$
  The models being identical,  we conclude that the algebras $A\*\Cr\B$ and
$\Cr{A\*\B}$ are naturally  isomorphic.
\endProof

\state Proposition  \label ExactReducedCrossSec
  Let $\B = \{B_g\}_{g\in G}$ be a Fell bundle over an exact discrete group $G$.  If
$B_1$ is an exact C*-algebra, then so is $\Cr\B$.

\Proof
  Let
  $$
  0 \to J \lto \iota A \lto \pi Q \to 0
  $$
  be an exact sequence of C*-algebras.  We need to prove that
  $$
  0 \to J\*\Cr\B \Lto{\iota\*1} A\*\Cr\B \Lto{\pi\*1} Q\*\Cr\B \to 0
  $$
  is also exact.
  Employing the isomorphisms obtained in \lcite{\TensorProdFBun} we may instead
prove the exactness of the sequence
  $$
  0 \to \Cr{J\*\B} \Lto{\iota\*1} \Cr{A\*\B} \Lto{\pi\*1} \Cr{Q\*\B} \to 0.
  \subeqmark TheExactSeq
  $$
  In naming the arrows in the above sequence we have committed a slight abuse of
language since we should actually have  employed the isomorphisms obtained in
\lcite{\TensorProdFBun}. Nevertheless, if the map we labeled $\pi\*1$ in the
last sequence above is applied to an element in $\Cr{A\*\B}$ of the form
$a\*b_g$, with $b_g\in B_g$, the result will be $\pi(a)\*b_g$, so we feel our choice of notation is justified.

  As it is well known,  the only possibly controversial point relating to the
exactness of \lcite{\TheExactSeq} is whether or not the kernel of
$\pi\*1$,  which we will refer to as $K$,  is contained in the image
of $\iota\*1$.
  We will arrive at this conclusion by applying \cite[5.3]{\exact} to $K$.  For
this we need to recall from \cite[3.5]{\amena} that, for each $g$ in $G$, there
is a contractive linear map
  $$
  F_g: \Cr\B \to B_g
  $$
  satisfying $F_g(\sum_{h}b_h) = b_g$, whenever $(b_h)_h$ is a
finitely supported section of $\B$.   Here we shall make use of these maps both for
the Fell bundle $A\*\B$ and for $Q\*\B$,  and  we will denote them by $F_g^A$ and
$F_g^Q$, respectively.

According to \cite[5.2]{\exact}, to check that the ideal $K$ in $\Cr{A\*\B}$ is
\"{invariant} we must verify that $F_g^A(K)\subseteq K$, for each $g$ in $G$.
For this we consider the diagram
  $$
  \matrix{
   \Cr{A\*\B} & \HighLto{\pi\*1} &    \Cr{Q\*\B} \cr
   \pilar{20pt}
   F_g^A\Big\downarrow\quad && \quad\Big\downarrow F_g^Q \cr
   \pilar{20pt}
   A\*B_g & \HighLto{\pi\*1}  &    Q\*B_g
   }
  $$
  In order to check that this is commutative,  let $x\in \Cr{A\*\B}$ have the
form
  $
  x = a\*b_h,
  $
  where $a\in A$ and $b_h\in B_h$, for some $h\in G$.  Employing Kronecker
symbols we then have that
  $$
  (\pi\*1)F^A_g(x) = \delta_{gh}(\pi\*1)(x) = \delta_{gh}\pi(a)\*b_h,
  $$
  while
  $$
  F^Q_g(\pi\*1)(x) = F^Q_g(\pi(a)\*b_h) = \delta_{gh}\pi(a)\*b_h.
  $$
  Since the set of elements $x$ considered above clearly generates $\Cr{A\*\B}$,
we see that the diagram is indeed commutative.
  If we now take an arbitrary element $x\in K$,  we will have that
  $$
  0 =   F^Q_g(\pi\*1)(x) = (\pi\*1)F^A_g(x),
  $$
  which implies that $F^A_g(x)\in K$, meaning that $K$ is invariant under $F^A_g$.

Given that $G$ is assumed to be exact,  we may apply \cite[5.3]{\exact} to conclude
that $K$ is \"{induced}, meaning that it
is generated, as an ideal, by its intersection with the unit fiber algebra,
namely $K\cap(A\*B_1)$.  The latter evidently coincides with the kernel of the
restriction of $\pi\*1$ to $A\*B_1$.  However, since the image of $A\*B_1$ under
$\pi\*1$ is contained in $Q\*B_1$, we may view $K\cap(A\*B_1)$ as the kernel of
the third map in the sequence
  $$
  0 \to J\*B_1 \Lto{\iota\*1} A\*B_1 \Lto{\pi\*1} Q\*B_1 \to 0.
  $$
  At this point we invoke our second main hypothesis, namely that $B_1$ is
exact, to deduce that the sequence above is exact, and hence that
$K\cap(A\*B_1)=J\*B_1$.
  Using angle brackets to denote generated ideals we then have that
  $$
  K =
  \big\langle K\cap(A\*B_1) \big\rangle =
  \big\langle J\*B_1\big\rangle \subseteq \Cr{J\*\B},
  $$
  which proves that
\lcite{\TheExactSeq} is exact in the middle. \endProof

\state Corollary Given a partial action $\alpha$ of an exact discrete group $G$ on an exact
C*-algebra $A$, the reduced crossed product $A\xr_{\alpha}G$ is
exact.

\Proof
  It is enough to notice that $A\xr_{\alpha}G$ is the reduced
cross-sectional C*-algebra of the semidirect product bundle, which is a Fell
bundle over $G$, and has $A$ as the unit fiber algebra.
  \endProof

Recalling from \lcite{\DefineRedOnm} that $\Omn^\red$ is the reduced crossed
product of an abelian,  hence exact,  C*-algebra by the exact free group $\Fmn$, we obtain:

\state Corollary For every positive integers $n$ and $m$,  one has that $\Omn^\red$
is an exact C*-algebra.

\section On full cross-sectional C*-algebras of Fell bundles

We shall now prove some preparatory results in order to study $\Omn$ (rather
than  the reduced version $\Omn^\red$).  Our goal is to show that it is not an
exact C*-algebra,  for $m, n\geq 2$,  from which it will follow that it indeed
differs from its reduced counterpart.

Since $\Omn$ is the \"{full} crossed product $\MCP$, we will now concentrate on
\"{full} cross-sectional algebras of Fell bundles.  However we will start with a
result about \"{reduced} cross-sectional algebras which will prove to be quite
useful in the study of their full versions.

\state Proposition \label SubFellBundes
  Let
  $\B = \{B_g\}_{g\in G}$
  be a Fell bundle over a discrete group $G$ and let $H$ be a subgroup of $G$.  Denote by
  $\C=\{C_h\}_{h\in H}$
  the Fell bundle obtained by restricting $\B$ to $H$, meaning that $C_h = B_h$,
for each $h\in H$, with norm, multiplication and involution borrowed from $\B$.
  Then: \izitem
  \zitem There exists a conditional expectation $E$ on $\Cr\B$ whose range is
isomorphic to $\Cr\C$.
  \zitem If $\Cr\B$ is nuclear (resp.~exact), then so is $\Cr\C$.

\Proof
  Viewing each $B_g$ as a subset of $\Cr\B$, as allowed by \cite[2.5]{\amena},
let $A$ be the closed linear span of $\bigcup_{h\in H} C_h$.  The
standard conditional expectation $E:\Cr\B \to B_1$ given by \cite[2.9]{\amena}
may be restricted to give a conditional expectation from $A$ to $C_1=B_1$,
satisfying the hypothesis of \cite[3.3]{\amena}.  Consequently there exists a
surjective *-homomorphism
  $$
  \lambda : A \to \Cr\C.
  $$
  By \cite[3.6]{\amena} the kernel of $\lambda$ is the set formed by the
elements $a\in A$ such that $E(a^*a)=0$.  However, applying \cite[2.12]{\amena}
to $\Cr\B$, one sees that only the zero element satisfies such an equation, which means that
$\lambda$ is injective and hence that $A$ is isomorphic to $\Cr\C$.

We now claim that the map
  $$
  E_H\ :\ \sum_{g\in G} b_g\ \in\  \bigoplus_{g\in G} B_g\
  \mapsto \ \sum_{g\in H} b_g\ \in \ A
  $$
  is continuous relative to the norm on its domain induced by $\Cr\B$.  In order
to see this recall that,
strictly according to definition \cite[2.3]{\amena}, $\Cr\B$ is the closed
*-subalgebra of
  $\Lin\B$
  (adjointable operators on the right Hilbert $B_1$--module $\ell_2(\B)$)
  generated by the range of the left regular representation of
$\B$.

Let $\iota$ be the natural inclusion of $\ell_2(\C)$ into $\ell_2(\B)$ and
observe that its adjoint is the projection of the latter onto the former.
Now consider the linear map
  $$
  V: T\in \Lin\B \mapsto \iota^* T \iota\in \Lin C.
  $$
  Viewing each $B_g$ within $\Lin\B$, and each $C_h$ within $\Lin\C$, by
\cite[2.2 \& 2.5]{\amena}, one may easily show that for every $g\in G$, and
every $b_g\in B_g$, one has that
  $$
  V(b_g) = \left\{
  \matrix{ b_g, & \hbox{if } g\in H,\hfill \cr \pilar{14pt} 0, & \hbox{otherwise.}
  }\right.
  $$
  Therefore, given any $\sum_{g\in G} b_g \in  \bigoplus_{g\in G} B_g$, we have that
  $$
  \Big\| \sum_{h\in H} b_h\Big\| =
  \Big\|V\Big(\sum_{g\in G} b_g\Big)\Big\| \leq
  \Big\|\sum_{g\in G} b_g\Big\|,
  $$
  where the norm in the left hand side is computed in $\Lin\C$.  However, due to
the fact that $A$ and $\Cr\C$ are isomorphic, the inequality above also holds
if the norm in the left hand side is computed in $A$.  This says that $E_H$ is
continuous, hence proving our claim.

One may then easily prove that the unique continuous extension of $E_H$ to $A$
is a conditional expectation, taking care of point (i).

Point (ii) now follows immediately from (i), since a closed *-subalgebra of an
exact C*-algebra is exact \cite[Exercise 2.3.2]{\Nate}, and the range of a conditional
expectation on a nuclear C*-algebra is nuclear \cite[Exercise 2.3.1]{\Nate}.
\endProof

Recall that $\l$ denotes
the left regular representation of a group $G$ in $\Cr G$.  Also, given a Fell bundle $\B$, we
will let $\Lambda$ be the regular representation of $\B$ in $\Cr\B$
\cite[2.2]{\amena}.

Recall from \cite[VIII.16.12]{\FD} that every representation $\pi$ of $\B$ in a
C*-algebra $A$ extends to a *-homomorphism (also denoted $\pi$ by abuse of
language) from $\Cf\B$ to $A$.

We thank Eberhard Kirchberg for sharing with us a very interesting idea which, when applied to
Fell bundles, yields the following curious result, mixing \"{reduced} cross-sectional
C*-algebras and \"{maximal} tensor products to produce \"{full} cross-sectional
C*-algebras.  See also \cite[10.2.8]{\Nate}.

\state Theorem \label CrazyFact
  Let $\Lambda\tm\l$ be the representation of\/ $\B$ in $\Cr\B \tm \Cr G$ given
by
  $$
  (\Lambda\tm\l)b_g = \Lambda(b_g) \* \l_g \for g\in G\for b_g\in B_g.
  $$
  Then the associated *-homomorphism
  $$
  \Lambda\tm\l: \Cf\B \to \Cr\B \tm \Cr G
  $$
  is injective.

\Proof Choose a faithful representation $\pi: \Cf\B \to B(H)$, where
$H$ is a Hilbert space, and consider the representation
$\pi\*\lambda$ of $\B$ on $H\*\ell_2(G)$ given by
  $$
  (\pi\*\l)b_g = \pi(b_g)\*\l_g \for g\in G\for b_g\in B_g.
  $$
  This gives
rise to the representation
  $\pi\*\lambda$ of $\Cf\B$
  which factors through a representation
  $$
  \pi_\lambda:\Cr\B \to B\big(H\*\ell_2(G)\big),
  $$
  by \cite[3.4]{\exact}.
  Let $\rho$ be the \"{right} regular representation of $G$ on $\ell_2(G)$,
which in turn yields the representation $\tilde \rho$ of $\Cr G$ on
$H\*\ell_2(G)$ defined by
  $$
  \tilde\rho = 1\*\rho:\Cr G \to B\big(H\*\ell_2(G)\big).
  $$
  It is easy to see that   the range of $\pi_\lambda$ commutes with the range of $\tilde\rho$, so there exists
a representation
  $$
  \pi_\lambda\*\tilde\rho : \Cr\B \tm\Cr G \to B\big(H\*\ell_2(G)\big),
  $$
  such that
  $$
  (\pi_\lambda\*\tilde\rho)(x\*y) = \pi_\lambda(x)\tilde\rho(y)
  \for x\in \Cr\B \for y\in \Cr G.
  $$
  Given any $g$ in $G$, and any $b_g\in B_g$, observe that
  $$
  (\pi_\lambda\*\tilde\rho)(\Lambda\tm\l)b_g =
  (\pi_\lambda\*\tilde\rho)\big(\Lambda(b_g)\*\l_g\big) =
  \pi_\lambda\big(\Lambda(b_g)\big) \ \tilde\rho(\l_g) \$=
  \big(\pi(b_g)\*\l_g\big) (1\*\rho_g) =
  \pi(b_g)\*\l_g\rho_g.
  $$
  Denoting by $\{\d_g\}_{g\in G}$ the standard orthonormal basis of $\ell_2(G)$,
pick any $\xi\in H$, and observe that the above operator,  when applied to $\xi\*\d_1$
produces
  $$
  (\pi_\lambda\*\tilde\rho)(\Lambda\tm\l)b_g\calcat{\xi\*\d_1} =
  \big(\pi(b_g) \* \l_g\rho_g\big)(\xi\*\d_1) =
  \pi(b_g)\xi \* \d_1.
  $$
  By linearity,  density and continuity we conclude that
  $$
  (\pi_\lambda\*\tilde\rho)(\Lambda\tm\l)x\calcat{\xi\*\d_1} =
  \pi(x)\xi \* \d_1
  \for x\in \Cf\B.
  $$
  Therefore,  assuming that $(\Lambda\tm\l)x=0$, for some $x\in \Cf\B$, we
deduce that   $\pi(x)\xi=0$, for all $\xi\in H$, and hence that $\pi(x)=0$.
Since $\pi$ was supposed to be injective on $\Cf\B$, we deduce that $x=0$.
\endProof

The following is also based on an idea verbally communicated to us by  Kirchberg.

\state Corollary \label SubBundle
  Let $\B$ be a Fell bundle over a discrete group $G$ and let $H$ be a subgroup
of $G$.  Consider the Fell bundle
  $\C=\{C_h\}_{h\in H}$
  obtained by restricting $\B$ to $H$, meaning that $C_h = B_h$,
for each $h\in H$, with norm, multiplication and involution borrowed from $\B$.
  Then the natural map $\iota : \Cf\C \to \Cf\B$ is injective.

\Proof
  Recall from \lcite{\SubFellBundes.i} that there exists a conditional expectation from $\Cr\B$ onto
$\Cr\C$, as well as a conditional expectation from $\Cr G$ to $\Cr H$.
Therefore by \cite[3.6.6]{\Nate} one has that the natural maps below are
injective:
  $$
  \Cr\C \tm \Cr H \ \hookrightarrow \ \Cr\B\tm \Cr H \ \hookrightarrow \ \Cr\B \tm \Cr G.
  $$
  Consider the diagram
  $$
  \matrix{
  \Cr\C \tm \Cr H &  \hookrightarrow  & \Cr\B \tm \Cr G \cr
  \pilar{16pt}  \uparrow && \uparrow \cr
  \pilar{20pt}   \Cf\C & \buildrel \iota \over \longrightarrow& \Cf\B}
  $$
  where the vertical arrows are the versions of $\Lambda\tm\l$ for $\C$ and
$\B$, respectively.  By checking on elements $c_h\in C_h$, it is elementary to
prove that the above diagram commutes.  Since  all arrows, with the possible
exception of $\iota$,
are known to be injective, we deduce that $\iota$ is injective as well.
\endProof

The following is an interesting conclusion to be drawn from \lcite{\CrazyFact}.

\state Theorem \label NuclearAmenable
  Let $\B$ be a Fell bundle over the discrete group $G$.  If the reduced
cross-sectional C*-algebra $\Cr\B$ is nuclear, then the regular representation
  $$
  \Lambda: \Cf\B \to \Cr\B
  $$
  is an isomorphism.

\Proof
  Consider the commutative diagram
  $$
  \matrix{
  \Cf\B &\Lto{\Lambda\tm\l} & \Cr\B \tm \Cr G \cr
  \pilar {20pt} \Lambda \Big\downarrow && \quad \Big\downarrow q\cr
  \pilar {20pt} \Cr\B &\HighLto{id \*\l }& \Cr\B \* \Cr G \cr
  }
  $$
  where $q$ is the natural map from the maximal to the minimal tensor product.
Assuming that $\Cr\B$ is nuclear we have that $q$ is injective
\cite[3.6.12]{\Nate}, and hence $\Lambda$ is injective.  \endProof

\state Remark \rm According to \cite[4.1]{\amena}, the above result says that
$\B$ is an \"{amenable} Fell bundle.  However, as observed in the
very last paragraph of \cite{\ortho}, we do not know whether this implies the
\"{approximation property} for $\B$ \cite[4.4]{\amena}.  Nevertheless,  in view of
\cite[4.4.3]{\Nate},  it is perhaps reasonable to believe that the approximation
property could be deduced from the nuclearity of $\Cr\B$.

\section Isotropy groups for partial actions

Given a partial action
  $$
  \theta = \big\{\theta_g:X_{g\inv}\to X_g\big\}_{g\in G}
  $$
  of a discrete group $G$ on a locally compact Hausdorff topological space $X$, recall
that the \"{isotropy subgroup} for a given point $x\in X$ is the subgroup of $G$
defined by
  $$
  G^x = \big\{g\in G: x\in X_{g\inv},\ \theta_g(x)=x\big\}.
  $$

\state Proposition \label LotsOfConclusions
  Let $X$ be a Hausdorff locally compact topological space, let $G$ be a
discrete group, and let $\theta$ be a partial action of $G$ on $X$.
 Then:
  \izitem
  \zitem If the full crossed product $C_0(X)\rtimes_{\theta}G$ is exact, then
for every $x$ in $X$ for which $G^x$ is {\resfinGrp} {\rm
\cite[p.~96]{\Nate}}, one has that $G^x$ is amenable.
  \zitem
  If the reduced crossed product $C_0(X)\xr_{\theta}G$ is nuclear, then
the isotropy group of every point in $X$ is amenable.
  \zitem   If the reduced crossed product $C_0(X)\xr_{\theta}G$ is exact,
then the isotropy group of every point in $X$ is exact.

\Proof
  Given $x$ in $X$, consider the restriction of $\theta$ to $G^x$, thus
obtaining a partial action of $G^x$ on $X$.  Observing that the full crossed
product is defined to be the full cross-sectional C*-algebra of the associated
semidirect product Fell bundle, we deduce from \lcite{\SubBundle} that
$C_0(X)\rtimes_{\theta}G^x$ is isomorphic to a closed *-subalgebra of
$C_0(X)\rtimes_{\theta}G$.  By the assumption in (i) that the latter is exact,
we deduce that $C_0(X)\rtimes_{\theta}G^x$ is also exact \cite[Exercise
2.3.2]{\Nate}.

Consider
the *-homomorphism
  $$
  \pi:f\in C_0(X) \mapsto f(x)\cdot 1\in\Cf{G^x},
  $$
  as well as the universal representation of $G^x$
  $$
  u: G^x \to \Cf{G^x}.
  $$
  Viewing $\Cf{G^x}$ as an algebra of operators on some Hilbert space, it is
easy to check that $(\pi,u)$ is
a covariant representation
of the partial dynamical system $\big(C_0(X),G^x,\theta|_{G^x}\big)$,  in the
sense of \cite[1.2]{\topfree}.
  Therefore, by \cite[1.3]{\topfree} there exists a *-homomorphism
  $$
  \pi\times u: C_0(X)\rtimes_{\theta}G^x \to \Cf{G^x}
  $$
  such that
  $$
  (\pi\times u)(f\delta_h) = f(x)u_h,
  $$
  for all $h$ in $G^x$, and all $f$ in $C_0(X_h)$.  One moment of reflexion is
enough to convince ourselves that $\pi\times u$ is surjective and hence that
$\Cf{G^x}$ is exact by \cite[9.4.3]{\Nate}.

Under the assumption that $G^x$  is {\resfinGrp} we then deduce from
\cite[3.7.11]{\Nate} that $G^x$ is amenable, completing the proof of (i).

  We next consider the diagram
  $$
  \matrix{
   C_0(X)\rtimes_{\theta}G^x  & \HighLto{\pi\times u} &   \Cf{G^x} \cr
   \pilar{20pt} \Lambda\big\downarrow\quad && \quad\big\downarrow \Lambda^x \cr
   \pilar{20pt} C_0(X)\xr_{\theta}G^x  & \buildrel \phi \over {\cdots\cdots} &   \Cr{G^x} \cr
   \pilar{20pt} E\big\downarrow\quad && \quad\big\downarrow \tau \cr
   \pilar{20pt} C_0(X) & \HighLto{\chi^x}  &   \Cplx
   }
  \subeqmark  SixDiagram
  $$
  where
  $\Lambda$ is the \"{left regular representation} (see the paragraph following
\cite[2.3]{\amena}),
  $\Lambda^x$ is the  version of $\Lambda$ for the trivial one-dimensional Fell
bundle over $G^x$,
  $E$ is the standard conditional expectation \cite[2.9]{\amena},
  $\tau$ is the unique normalized trace on $\Cr{G^x}$ such that $\tau(\l_h) =
0$, for all $h\neq 1$, and finally
  $\chi^x$ is the character on $C_0(X)$ given by point evaluation at $x$.
  Incidentally $\tau$ coincides with the standard conditional expectation in
the context of the trivial bundle over $G^x$.

By checking on elements of the form
$f\delta_h$, it is elementary to verify that the diagram commutes.
We claim that $\pi\times u$ maps the kernel of $\Lambda$ into the kernel of
$\Lambda^x$.  In order to see this, suppose that $x$ lies in the kernel of
$\Lambda$.    Then by \cite[3.6]{\amena} we have that $E(x^*x)=0$, so
  $$
  0 = \chi^x\big(E\big(\Lambda(x^*x)\big)\big) =
  \tau\big(\Lambda^x\big((\pi\times u)(x^*x)\big)\big) =
  \tau\big(y^*y),
  $$
  where $y = \Lambda^x\big((\pi\times u)x\big)$.
  Since $\tau$ is a faithful trace on $C^*_{\red}(G^x)$ \cite[2.12]{\amena}, we
conclude that $y=0$, which proves that $(\pi\times u)x$ belongs to the kernel of
$\Lambda^x$, hence the claim.

  As a consequence we see that there exists a *-homomorphism $\psi$
  filling the dots in \lcite{\SixDiagram} in a way as to preserve the
commutativity of the diagram.  Since $\Lambda^x$ is surjective, $\psi$ must also
be surjective.

  Assuming that $C_0(X)\xr_{\theta}G^x$ is nuclear (resp.~exact), we
now deduce that $C^*_{\red}(G^x)$ shares this property.  To conclude the proof
it is now enough to recall that
  if $C^*_{\red}(G^x)$ is nuclear then $G^x$ is amenable \cite[2.6.8]{\Nate},
and that
  if $C^*_{\red}(G^x)$ is an exact C*-algebra then $G^x$ is an exact group \cite[5.1.1]{\Nate}.
\endProof

\state Theorem \label OnmNotExact
  If $m,n\ge2$, then $\Omn$ is not exact and hence it is not isomorphic to $\Omn^\red$.

\Proof
  Recall from \lcite{\FixedPoint} that there exists $y$ in $Y\univ$ such that
  $$
  (v\univ_1)\inv h\univ_1(y) = y = (v\univ_2)\inv h\univ_2(y).
  $$
  This implies that $b_1\inv a_1$ and $b_2\inv a_2$ belong to $\F^y_{m+n}$, the
isotropy group of $y$.

It is easy to see that the subgroup of $\Fmn$ generated by  these two elements
is isomorphic to $\F_2$, so we conclude that $\F^y_{m+n}$ is not amenable.

It is well known that free groups are {\resfinGrp} \cite[Corollary
22]{\Cohen} and consequently the same applies to its subgroup $\F^y_{m+n}$.
  Using \lcite{\LotsOfConclusions.i} one deduces that the full crossed product
$\MCP$ cannot be exact,  and hence the conclusion then follows from \lcite{\IsoOnmProd}.
   \endProof

% \---input /Users/exel/Casa/z/q/ara-katsura/findimrep.tex

\section Absence of finite dimensional representations

The goal of this section is to prove that $\Omnr$ does not admit any nonzero
finite dimensional representation.   In case $n\neq m$ the same is true
even for the unreduced algebras and, since the proof of this fact is much simpler, we
present it first.

\state Proposition \label EasyInexistFinDimRep
  If $n\neq m$ then $\Omn$ (and hence also $\Omnr$) does not admit any
nontrivial finite dimensional representation.

\def\uv{\underline v}
\def\uw{\underline w}
\def\tr{{\rm tr}}

\Proof
  Let $\rho:\Omn \to M_d(\Cplx)$ be a non-degenerate $d$-dimensional
representation,  with $d>0$.  Then,
denoting by $\uv$ and $\uw$ the images of $v$ and $w$ (see the third and fourth
relation in $(\R)$), we have
  $$
  \tr\big(\rho(\uv)\big)  = \sum_{i=1}^n \tr\big(\rho(\us_i\us_i^*)\big) =
  \sum_{i=1}^n \tr\big(\rho(\us_i^*\us_i)\big) = n\, \tr\big(\rho(\uw)\big),
  $$
  and similarly
  $
  \tr\big(\rho(\uv)\big)  = m\, \tr\big(\rho(\uw)\big),
  $
  so
  $$
  n\, \tr\big(\rho(\uw)\big) = m\, \tr\big(\rho(\uw)\big).
  $$
  Since $n\neq m$,  this implies that $\tr\big(\rho(\uw)\big)=0$, and hence also that
$\tr\big(\rho(\uv)\big)=0$.  Therefore
  $$
  d = \tr(1) = \tr(\rho(1)) = \tr\big(\rho(\uv+\uw))=0,
  $$
  a contradiction.
  \endProof

From now on we will develop a series of auxiliary results in order to
show the nonexistence of nonzero finite dimensional representations of $\Omnr$
when $m=n$ (although our proof will not explicitly use that $m=n$,
and hence it will serve as a proof for the general case).
In what follows we will therefore assume that
  $$
  \rho:\Omnr \to M_d(\Cplx)
  $$
  is non-degenerate $d$-dimensional representation and our task will be to
arrive at a contradiction from it.

  Restricting $\rho$ to $C(\OR)$ we get a finite dimensional representation of a
commutative algebra which, as it is well known, is equivalent to a direct sum of
characters.  In other words, upon conjugating $\rho$ by some unitary matrix, we
may assume that there is a $d$-tuple $(\xi_1,\xi_2,\ldots,\xi_d)$ of elements of
$\OR$ such that
  $$
  \rho(f) = \pmatrix{ f(\xi_1) \cr & f(\xi_2) \cr && \ddots \cr &&& f(\xi_d)},
  $$
  for every $f$ in $C(\OR)$.

\state Proposition The set $Z =\{\xi_1,\xi_2,\ldots,\xi_d\}$ is invariant under
$\tu$.

\Proof We want to prove that for every $g$ in $\Fmn$, and every
  $\xi\in Z\cap \OR_{g\inv}$, one has that $\tu_g(\xi)$ is in $Z$.  Arguing by
contradiction we assume that this is not so, that is, that we can find
  $\xi\in Z\cap \OR_{g\inv}$ such that $\tu_g(\xi)\notin Z$.  Observing that $\tu_g(\xi)\in
\OR_g$, we may pick an $f\in C_0(\OR_g)$ such that $f(\tu_g(\xi))$ is nonzero,
but such that $f$ vanishes identically on $Z$.  In particular this implies that
$\rho(f)=0$.

Using \cite[1.4]{\topfree} we may write $\rho = \pi\times u$, where $(\pi,u)$ is
a covariant representation of the dynamical system $\big(C(\OR),\Fmn,\tu\big)$.
Noticing that $\pi$ is the restriction of $\rho$ to $C(\OR)$, we have
  $$
  \rho\big(\tu_{g\inv}(f)\big) =
  \pi\big(\tu_{g\inv}(f)\big) = u_{g\inv}\pi(f) u_g =0.
  $$
  It follows that
  $$
  0 = \tu_{g\inv}(f)\calcat {\xi} = f\big(\tu_g(\xi)\big)  \neq 0,
  $$
  a contradiction.
  \endProof

\state Proposition \label FreeIsotropy
  If $m, n\geq 2$, then for every $\xi$ in $Z$, the isotropy group
$\F^\xi_{m+n}$, contains a subgroup isomorphic to $\F_2$.

\Proof
  Assume first that $\xi\in Y\univ$, that is, the configuration of $\xi$ at
the origin follows pattern ($c_2$).  Then in particular $b_1\inv\in \xi$, and
hence the configuration of $\xi$ at $b_1$ must follow pattern ($c_1$).
Therefore there exists a unique $i_1\leq n$, such that $b_1\inv a_{i_1}\in \xi$.
The configuration of $\xi$ at $b_1\inv a_{i_1}$ must then follow pattern ($c_2$)
so, in particular $b_1\inv a_{i_1}b_1\inv\in\xi$.

Continuing in this way we may construct an infinite  sequence $i_1,i_2,\ldots$
such that
  $$
  g_k := b_1\inv a_{i_1} b_1\inv a_{i_2}b_1\inv  \ldots b_1\inv a_{i_k}\in\xi
  \for k\in\N.
  $$
  So $\xi\in \OR_{g_k}$, and hence
  $$
  \tu_{g_k\inv}(\xi) = g_k\inv\xi \in Z,
  $$
  because $Z$ is invariant under $\tu$.  Using the fact that $Z$ is finite we
conclude that there are positive integers $k<l$, such that
  $$
  g_l\inv\xi  = g_k\inv\xi,
  $$
  so $g_kg_l\inv\xi=\xi$, and hence the element
  $$
  x := g_kg_l\inv
  $$
  lies in the isotropy group of $\xi$.

Let $\F_2$ be the free group on a set of two generators, say $\{c_1,c_2\}$, and
consider the unique group homomorphism
  $$
  \phi: \Fmn \to \F_2
  $$
  such that
  $$
  \phi(a_i) = 1 \for i=1,\ldots,n,
  $$
  $$
  \phi(b_1)=c_1, \quad \phi(b_2)=c_2,
  \quad \phi(b_j) =1\for j\geq 3.
  $$

It is then evident that $\phi(g_k) = c_1^{-k}$, and hence that
  $$
  \phi(x) = \phi(g_kg_l\inv) = c_1^{l-k},
  $$
  where by assumption, $l-k>0$.

Repeating the above argument with $b_2$ in place of $b_1$, we may find some $y$
in the isotropy group of $\xi$ such that $\phi(y)$ is a positive power of $b_2$.

The subgroup of $\F^\xi_{m+n}$ generated by $x$ and $y$
is therefore a free group since its image within $\F_2$ via $\phi$ is certainly
free.

This concludes the proof under the assumption that the configuration of $\xi$ at
the origin is ($c_2$), so let us  suppose that the pattern is
($c_1$).
   Therefore there exists some $i$ such that $a_i\in\xi$ and hence,  again by invariance
of $Z$, we have that $a_i\inv\xi\in Z$.  Since $1\in\xi$ we have that $a_i\inv
\in a_i\inv\xi$, so the pattern of $a_i\inv\xi$ at the origin is necessarily ($c_2$).

By the case already studied there is a copy of $\F_2$ inside the isotropy group
of $a_i\inv\xi$, but since
  $$
  \F^{a_i\inv\xi}_{m+n} =   a_i\inv(\F^\xi_{m+n})a_i,
  $$
  the same holds for the isotropy group of $\xi$.
\endProof

Since $Z$ is invariant under $\tu$ we may restrict the latter to the former thus
obtaining a partial action,  say $\theta$, of $\Fmn$ on $Z$.
  % and so we may form the crossed product $C(Z)\xr_\theta\Fmn$.

Given $\xi\in Z$, we will denote by $1_\xi$ the characteristic function of
the singleton $\{\xi\}$, viewed as an element of $C(Z)$.

\state Proposition \label CopyCstarFtwo
  For every $\xi\in Z$ there exists an embedding of $\Cr{\F_2}$ in the reduced
crossed product $C(Z)\xr_\theta\Fmn$, such that the unit of the former is mapped
to $1_\xi$.

\Proof Let $G$ be any subgroup of $\F^\xi_{m+n}$.
For each $g$ in $G$, consider the element
  $$
  u_g =  1_\xi\d_g \in C(Z)\xr_\theta\Fmn.
  $$
  By direct computation one checks that
  $u_gu_h = u_{gh}$,  and
  $u_{g\inv} = u_g^*$,
  for every $g$ and $h$ in $G$, and moreover that
  $u_1 = 1_\xi$.  In other words, $u$ is a unitary representation of $G$ in the
hereditary subalgebra of $C(Z)\xr_\theta\Fmn$ generated by $1_\xi$.  Let
  $$
  \phi: C^*(G) \to C(Z)\xr_\theta\Fmn
  $$
  be the integrated form of $u$.  Denoting by $\tau$ the canonical trace on
$C^*(G)$, and by $E$ the standard conditional expectation
  $$
  E:C(Z)\xr_\theta\Fmn \to C(Z),
  $$
  one may easily prove that
  $$
  E\big(\phi(x)\big)= \tau(x)1_\xi  \for x\in C^*(G).
  $$
  Since $E$ is faithful, for every $x\in C^*(G)$ one has that
  $$
  \phi(x) = 0 \iff
  E(\phi(x^*x))=0 \iff
  \tau(x^*x)=0.
  $$
  This said we see that the kernel of $\phi$ coincides with the kernel of the
integrated form of the left regular representation, namely
  $$
  \lambda: C^*(G) \to \Cr{G}.
  $$
  Consequently $\phi$ factors through $\Cr{G}$, providing a *-homomor\-phism
  $$
  \tilde \phi: \Cr{G} \to C(Z)\xr_\theta\Fmn,
  $$
  which is injective because of the above equality of null spaces.  Clearly
$\tilde\phi(1)=1_\xi$, as stated.  To conclude the proof it is therefore enough
to choose $G$ to be the subgroup of $\F^\xi_{m+n}$ given by
\lcite{\FreeIsotropy}.
  \endProof

The next significant step in order to obtain a contradiction from the existence
of $\rho$ is to prove that it admits a factorization
  $$
  \matrix{\Omnr & {\buildrel {\textstyle\rho} \over {\hbox to 30pt{\rightarrowfill}}}& M_d(\Cplx) \cr
  \pilar{20pt}\hfill \phi\searrow\kern-12pt && \kern-12pt \nearrow \tilde\rho \hfill \cr
  \pilar{20pt}& C(Z)\xr_\theta\Fmn}
  \eqmark FactorizatRho
  $$
  such that $\phi(f) = f|_{Z}$, for all $f\in C(\OR)$.

The poof of this
factorization may perhaps be of independent interest, so we prove it in a more
general context in the next section.   Although it may not look a very deep
result we have not been able to prove it in full generality, since we need
to use the exactness of free groups.

\section Invariant ideals

Let $G$ be a discrete group and let $\a$ be a partial action of $G$ on $A$.
For each $g$ in $G$,  denote by $A_g$ the range of $\a_g$.

\definition
  A closed two-sided ideal $K \trianglelefteq A$ is said to be $\a$-invariant if
  $$
  \a_g(K\cap A_{g\inv}) \subseteq K
  \for g\in G.
  $$
  Given such an ideal, let $B = A/K$, and denote the quotient map by
  $$
  q: A \to B.
  $$
  For each $g$ in $G$, consider the closed two-sided ideal of $B$ given by
$B_g=q(A_g)$. Given any $b\in B_{g\inv}$, write $b = q(a)$, for some $a\in
A_{g\inv}$, and define
  $$
  \b_g(a) := q\big(\a_g(a)\big).
  $$
  It is then easy to see that $\b_g$ becomes a *-isomorphism from $B_{g\inv}$ to
$B_g$,  also known as a \"{partial automorphism} of $B$.

\state Proposition \label QuotientPartialAction
  The collection of partial automorphisms $\{\b_g\}_{g\in G}$ forms a partial
action of $G$ on $B$.

\Proof
  If $I$ and $J$ are closed two-sided ideals of $A$, it is well known that every
element $z\in I\cap J$ may be written as a product $z=xy$, with $x\in I$, and
$y\in J$.  In other words
  $I\cap J = IJ$.  Therefore
  $$
  q(I\cap J) = q(IJ) = q(I)q(J) = q(I)\cap q(J).
  $$
  We then conclude that
  $$
  \b_g(B_{g\inv} \cap B_h) =
  \b_g\big(q(A_{g\inv}) \cap q(A_h)\big) =
  \b_g\big(q(A_{g\inv}\cap A_h)\big) =
  q\big(\a_g(A_{g\inv}\cap A_h)\big) \$=
  q\big(A_g\cap A_{gh}\big) =
  q(A_g)\cap q(A_{gh}) =
  B_g\cap B_{gh}.
  $$
  We leave the verification of the remaining axioms
  (\cite{\ExelCircle},
  \cite{\McClanahan},
  \cite{\tpa}) to the reader.
\endProof

\state Proposition \label CovariantQuotientMap
  Under the above assumptions, there exists a unique surjective *-homomorphism
  $$
  \phi: A\xr_\a G \to B\xr_\b G,
  $$
  such that $\phi(a_g\d_g) = q(a_g)\d_g$, for all $g\in G$, and all $a_g\in A_g$.

\font \rssmall = rsfs10 scaled 600
\def\A{\curly{A}}
\def\B{\curly{B}}

\def\As{{\hbox{\rssmall A}}}
\def\Bs{{\hbox{\rssmall B}}}

\Proof
  Recalling that the reduced crossed product C*-algebra coincides with the
reduced cross-sectional C*-algebra of the corresponding semidirect product
bundle \cite[2.8]{\tpa}, denote by $\A$ and $\B$ the corresponding Fell
bundles.
Precisely $\A = \big\{A_g\d_g\big\}_{g\in G}$, with multiplication
  $$
  (a_g\d_g,b_h\d_h) \in A_g\d_g \times A_h\d_h \mapsto
  \a_g\big(\a_g\inv(a_g)b_h\big)\d_{gh}\in A_{gh}\d_{gh}
  $$
  and involution
  $$
  a_g\d_g\in A_g\d_g  \mapsto \a_{g\inv}(a_g^*)\d_{g\inv} \in A_{g\inv}\d_{g\inv},
  $$
  and likewise for $\cal B$.
  It is then easy to see that the correspondence
  $$
  a_g\d_g \in A_g\d_g \mapsto q(a_g)\d_g \in B_g\d_g
  $$
  defines a homomorphism in the category of Fell bundles and hence induces a
*-homo\-morphism  of full cross-sectional C*-algebras
  $$
  \psi : C^*(\A) \to C^*(\B).
  $$
  Denoting by
  $$
  E: C^*(\A) \to A \and F:C^*(\B) \to B
  $$
  the corresponding conditional expectations \cite[2.9]{\amena}, one easily verifies that
$F\psi = q E$.  From this it follows that, for every element $x$ in the kernel of the
regular representation \cite[2.2]{\amena},
  $$
  \Lambda_\As : C^*(\A) \to \Cr{\A},
  $$
  one has that
  $$
  F\big(\psi(x^*x)\big) =  q\big(E(x^*x)\big) = 0,
  $$
  by \cite[3.6]{\amena}.  Therefore, by \cite[2.12]{\amena}, we see that
$\psi(x^*x)$ lies in the kernel of the regular representation $\Lambda_\Bs$
relative to $\cal B$.  We conclude that $\psi$ factors through the quotient
providing a map $\phi$ such that the diagram below is commutative.
  $$
  \matrix{
   C^*(\A) & \buildrel {\textstyle \psi} \over \longrightarrow & C^*(\B) \cr
   \pilar{20pt}
   \Lambda_\As\Big\downarrow\quad && \quad\Big\downarrow \Lambda_\Bs  \cr
   \pilar{20pt}
   \Cr{\A} & \buildrel {\textstyle \phi} \over \longrightarrow & \Cr{\B}.}
  $$
  Identifying reduced crossed products with their corresponding reduced
cross-sectional algebras, the proof is complete.
  \endProof

\state Proposition \label RepsFactor
  Let $\a$ be a partial action of a discrete {exact} group $G$ on a C*-algebra
$A$, and let $\rho$ be a *-representation of $A\xr_\a G$ on a Hilbert space $H$.
Letting $K$ be the null-space of $\rho|_A$, then $K$ is $\a$-invariant, so we
may speak of the \"{quotient partial action} $\b$ of
\lcite{\QuotientPartialAction}, and of the map $\phi$ of
\lcite{\CovariantQuotientMap}.  Under these conditions there exists a
*-representation $\tilde\rho$ of $A/K\xr_\b G$, such that the diagram
  $$
  \matrix{A\xr_\a G & {\buildrel {\textstyle\rho} \over {\hbox to 30pt{\rightarrowfill}}}& \B(H) \cr
  \pilar{20pt}\hfill \phi\searrow\kern-12pt && \kern-12pt \nearrow \tilde\rho \hfill \cr
  \pilar{20pt}& A/K\xr_\b G}
  $$
  commutes.

\Proof Let $J$ be the null space of $\rho$, so that $K=A\cap J$.  Given any
$g\in G$, and any
$a\in K\cap A_{g\inv}$, observe that, identifying $A$ with its image
in $A\xr_\a G$, as usual, one has that
  $$
  (b\d_g) a(b\d_g)^* = b\a_g(a)b^* \for b\in A_g.
  $$
  Applying $\rho$ on both sides of the above equality, we conclude that
$b\a_g(a)b^*\in K$.  If we now let $b$ run along an approximate identity for
$A_g$, we conclude that $\a_g(a)$ lies in $K$, thus proving that $K$ is
$\a$-invariant.

  We next claim that
  $$
  \hbox{Ker}(\phi) \subseteq \hbox{Ker}(\rho).
  \subeqmark KerContained
  $$
  With that goal in mind, let
  $$
  E: A\xr_\a G \to A
  \and
  F: A/K\xr_\b G \to A/K,
  $$
  be the associated conditional expectations (unlike
\lcite{\CovariantQuotientMap}, here these are seen as maps on the \"{reduced}
cross-sectional algebras).
  Given $x$ in the kernel of $\phi$, we have that
  $$
  0 =  F\big(\phi(x^*x)\big) = q\big(E(x^*x)\big),
  $$
  so we see that $E(x^*x)$ lies in $K \subseteq J$.  Using \cite[5.1]{\exact} we
deduce that $x$ is in the ideal of $A\xr_\a G$ generated by $K$, and hence that
$x$ is in $J$.  This proves \lcite{\KerContained} and, since $\phi$ is
surjective, we have that $\rho$ factors through $\phi$, which means precisely
that a map $\tilde\rho$ exists with the stated properties.
  \endProof

Returning to the situation we left at the end of the previous section, recall
that $\rho$ is a non-degenerate $d$-dimensional representation of $\Omnr$.
Notice that
  $$
  K:= \hbox{Ker}\big(\rho|_{C(\OR)}\big) =
  \big\{f\in C(\OR) : f(\xi_i)=0,\ \forall i=1,\ldots d\big\}.
  $$
  The quotient of $C(\OR)$ by $K$ may then be naturally identified with $C(Z)$,
and the {quotient partial action} given by \lcite{\QuotientPartialAction}
becomes the action induced by the restriction of $\tu$ to $Z$.  Thus, when
applied to our situation, the diagram in the statement of \lcite{\RepsFactor}
becomes precisely \lcite{\FactorizatRho}.

   The restriction of $\tilde\rho$ to the copy of
$\Cr{\F_2}$ provided by \lcite{\CopyCstarFtwo} will then be a (possibly degenerate)
$d$-dimensional representation of the \"{simple infinite-dimensional} C*-algebra
$\Cr{\F_2}$.  Such a representation must therefore be identically zero and hence, in
particular,
  $$
  \tilde\rho(1_\xi)=0,
  $$
  because,  as seen above,  $1_\xi$ lies in
the copy of $\Cr{\F_2}$ alluded to.
  Observing that the unit of $C(Z)\xr_\theta\Fmn$ is given by
  $$
  1 = \sum_{\xi\in Z}1_\xi,
  $$
  we deduce that $\tilde\rho(1)=0$ and hence that $\tilde\rho=0$.  A glance at
\lcite{\FactorizatRho} then gives $\rho=0$.

  This proves the following main result:

\state Theorem \label NoFinDimRepReg
  $\Omnr$ admits no nonzero finite dimensional representations.

  \endgroup

% \vskip 20pt

\Address
  {Departament de Matem\`atiques, Universitat Aut\`onoma de Barcelona, 08193
Bellaterra (Barcelona), Spain.}
  {para@mat.uab.cat}

\Address
  {Departamento de Matem\'atica, Universidade Federal de Santa Catarina,
88010-970 Florian\'opolis SC, Brazil.}
  {exel@mtm.ufsc.br}

\Address
  {Department of Mathematics, Faculty of Science and Technology, Keio
University, 3-14-1 Hiyoshi, Kouhoku-ku, Yokohama, Japan, 223-8522.}
  {katsura@math.keio.ac.jp}


\references

\bibitem{\AranAbrams}
  {G. Aranda Pino and G. Abrams}
  {The Leavitt path algebra of a graph}
  {\it J. Algebra \bf 293 \rm (2005),  319--334}

\bibitem{\AraGood1}
  {P. Ara and  K. R. Goodearl}
  {Leavitt path algebras of separated graphs}
  {to appear in J. reine angew. Math.;     arXiv:1004.4979v1 [math.RA] (2010)}


\bibitem{\AraGood}
  {P. Ara and  K. R. Goodearl}
  {C*-algebras of separated graphs}
  {\it J. Funct. Analysis \bf 261 \rm (2011), 2540--2568}

\bibitem{\ArMorPa}
   {P. Ara, M. A. Moreno, E. Pardo}
   {Nonstable $K$-theory for graph algebras}
   {\it Algebr. Represent. Theory \bf 10 \rm (2007), 157--178}

\bibitem{\LBrown}
 {L. G. Brown}
 {Ext of certain free product C*-algebras}
  {\it J. Operator Theory \bf 6 \rm (1981), 135--141}


\bibitem{\Nate}
  {N. P. Brown and N. Ozawa}
  {C*-algebras and finite-dimensional approximations}
  {Graduate Studies in Mathematics, 88, American Mathematical Society, 2008}

\bibitem{\Cohen}
  {D. E. Cohen}
  {Combinatorial group theory: a topological approach}
  {London Mathematical Society Student Texts, 14. Cambridge University Press, 1989}

\bibitem{\CuntzKrieger}
  {J. Cuntz and  W. Krieger}
   {A class of C*-algebras and topological Markov chains}
   {\it Invent. Math.  \bf 56  \rm (1980), 251--268}


\bibitem{\Pat}
  {J. Duncan and  A. L. T. Paterson}
  {C*-algebras of inverse semigroups}
  {\it Proc. Edinburgh Math. Soc. (2) \bf 28 \rm (1985), no. 1, 41--58}

\bibitem{\ExelCircle}
  {R. Exel}
  {Circle actions on C*-algebras, partial automorphisms and a generalized
Pimsner--Voiculescu exact sequence}
  {\it J. Funct. Analysis, \bf 122 \rm (1994), 361--401}

\bibitem{\amena} % \def\amena{E1}
  {R. Exel}
  {Amenability for Fell bundles}
  {\sl J. Reine Angew. Math. \bf 492 \rm (1997), 41--73
[arXiv:funct-an/9604009]}

\bibitem{\ortho} % \def\ortho{E2}
  {R. Exel}
  {Partial representations and amenable Fell bundles over free groups}
  {\it Pacific J. Math. \bf 192 \rm  (2000), 39--63
[arXiv:funct-an/9706001]}

  \bibitem{\exact} % \def\exact{E3}
  {R.~Exel}
  {Exact groups, induced ideals, and Fell bundles\fn{This paper was published as
``Exact groups and Fell bundles", \it Math.~Ann.~\bf 323 \rm (2002), no.~2,
259--266. However, the referee required that the results pertaining to induced
ideals be removed from the preprint version arguing that there were no applications of this concept.  The
reader will therefore have to consult the arxiv version, where the results we
need may be found.}}
  {\hfill\break http://arxiv.org/abs/math.OA/0012091 (2003)}

  \bibitem{\infinoa} % \def\infinoa{EL}
  {R. Exel and M. Laca}
  {Cuntz--Krieger algebras for infinite matrices}
  {\it J. reine angew. Math. \bf 512 \rm (1999), 119--172
[arXiv:funct-an/9712008]}

  \bibitem{\topfree} % \def\topfree{ELQ}
  {R.~Exel,  M.~Laca and John Quigg}
  {Partial dynamical systems and C*-algebras generated by partial
isometries\fn{This paper was published in
\it J.~Operator Theory \bf 47 \rm (2002), no.~1, 169--186.
However, the referee required that the results pertaining to covariant
representations be removed from the preprint version.  The
reader will therefore have to consult the arxiv version, where the results we
need may be found.}}
  {http://arxiv.org/abs/funct-an/9712007v1 (1997)}

\bibitem{\FD}
  {J. M. G. Fell and R. S.  Doran}
  {Representations of *-algebras, locally compact groups, and Banach *-algebraic
bundles. Vols. 1 \& 2}
  {Pure and Applied Mathematics, 126, Academic Press, 1988}

\bibitem{\Hazrat}
   {R. Hazrat}
   {The graded structure of Leavitt path algebras}
   {arXiv:1005.1900v4 [math.RA]}

\bibitem{\KR}
  {E. Kirchberg, M. R\o rdam}
  {Infinite non-simple $C^*$-algebras: Absorving the Cuntz algebra ${\cal
O}_{\infty}$}
  {\it Advances Math. \bf 167 \rm (2002), 195--264}

\bibitem{\Leavitt}
  {W. G. Leavitt}
  {The module type of a ring}
  {\it Trans. Amer. Math. Soc. \bf 103 \rm (1962), 113--130}

\bibitem{\McClanahanUnitary}
 {K. McClanahan}
 {C*-algebras generated by elements of a unitary
 matrix}
 {\it  J. Funct. Anal. \bf 107 \rm  (1992), 439–-457}

\bibitem{\McClanahanReactangle}
 {K. McClanahan} {K-theory and Ext-theory for rectangular unitary
 C*-algebras}
  {\it Rocky Mountain J. Math. \bf 23 \rm (1993), 1063--1080}

\bibitem{\McClanahanAmalgamated}
   {K. McClanahan}
   {Simplicity of reduced amalgamated products of C*-algebras}
   {\it Canad. J. Math. \bf 46 \rm  (1994), 793--807}

\bibitem{\McClanahan}
  {K. McClanahan}
  {$K$-theory for partial crossed products by discrete
groups}
  {\it J. Funct. Analysis \bf 130 \rm (1995), 77--117}
\bye